\documentclass[reqno]{amsart}
\usepackage{amssymb,wasysym,mathscinet,enumitem,bm}

\usepackage{mathrsfs}   
 \numberwithin{equation}{section}
 
\usepackage{geometry}
\usepackage[numbers,sort&compress]{natbib}
\usepackage{hyperref} 

\usepackage{esint} 

\theoremstyle{plain}
\newtheorem{theorem}{Theorem}[section]
\newtheorem{lemma}[theorem]{Lemma}

\newtheorem{proposition}[theorem]{Proposition}

\theoremstyle{definition}

\newtheorem{assumption}[theorem]{Assumption}

\theoremstyle{remark}
\newtheorem{remark}[theorem]{Remark}

 \makeatletter
 \newcommand{\norm}{\@ifstar{\@normb}{\@normi}}
 \newcommand{\@normb}[2]{\left\Vert{#1}\right\Vert_{#2}}
 \newcommand{\@normi}[2]{\Vert{#1}\Vert_{#2}}
 \makeatother

 \global\long\def\Sob#1#2{{W}^{#1}_{#2}} 
 \global\long\def\fSob#1#2{{H}^{#1}_{#2}}

 \global\long\def\oSob#1#2{\mathring{W}^{#1}_{#2}}
 \global\long\def\Leb#1{L_{#1}}

   \DeclareMathOperator{\Div}{div}
 \newcommand{\relphantom}[1]{\mathrel{\phantom{#1}}}
 
 \newcommand{\myd}[1]{\,d{#1}}

 \DeclareMathOperator{\supp}{supp}
 
 \DeclareMathOperator{\Tr}{Tr}

\makeatletter
\def\@tocline#1#2#3#4#5#6#7{\relax
  \ifnum #1>\c@tocdepth 
  \else
    \par \addpenalty\@secpenalty\addvspace{#2}%
    \begingroup \hyphenpenalty\@M
    \@ifempty{#4}{%
      \@tempdima\csname r@tocindent\number#1\endcsname\relax
    }{%
      \@tempdima#4\relax
    }%
    \parindent\z@ \leftskip#3\relax \advance\leftskip\@tempdima\relax
    \rightskip\@pnumwidth plus4em \parfillskip-\@pnumwidth
    #5\leavevmode\hskip-\@tempdima
      \ifcase #1
       \or\or \hskip 1em \or \hskip 2em \else \hskip 3em \fi%
      #6\nobreak\relax
    \hfill\hbox to\@pnumwidth{\@tocpagenum{#7}}\par
    \nobreak
    \endgroup
  \fi}
\makeatother

\usepackage{tcolorbox}\tcbuselibrary{breakable}

\allowdisplaybreaks

\begin{document}

\title[Stokes equations]{Nonstationary Stokes equations on a domain with curved boundary under slip boundary conditions}

\author[H. Dong]{Hongjie Dong}
\address{Division of Applied Mathematics, Brown University, 182 George Street, Providence, RI 02912, USA}
\email{hongjie\_dong@brown.edu }

\author[H. Kwon]{Hyunwoo Kwon}
\address{Division of Applied Mathematics, Brown University, 182 George Street, Providence, RI 02912, USA}
\email{hyunwoo\_kwon@brown.edu }
\thanks{H. Dong and H. Kwon were partially supported by the NSF under agreement DMS-2350129.}


\begin{abstract}
We consider nonstationary Stokes equations in nondivergence form with variable viscosity coefficients and generalized Navier slip boundary conditions with slip tensor $\mathcal{A}$ in a domain $\Omega$ in $\mathbb{R}^d$. First, under the assumption that slip matrix $\mathcal{A}$ is sufficiently smooth, we establish a priori local regularity estimates for solutions near a curved portion of the domain boundary. Second, when $\mathcal{A}$ is the shape operator, we derive local boundary estimates for the Hessians of the solutions, where the right-hand side does not involve the pressure. Notably, our results are new even if the viscosity coefficients are constant.
\end{abstract}

\maketitle

\section{Introduction}

This paper concerns local regularity estimates of solutions to nonstationary Stokes equations near a curved portion of the boundary of a domain. Let us consider the following nonstationary Stokes equations with variable viscosity coefficients:
\begin{equation}\label{eq:Stokes}
\left\{\begin{alignedat}{2}
\partial_t u -a^{ij}D_{ij} u +\nabla p &=f&&\qquad \text{in } (S,T)\times\Omega,\\
\Div u &=g&&\qquad \text{in } (S,T)\times\Omega.
\end{alignedat}
\right.
\end{equation}
Here $\Omega$ is a domain in $\mathbb{R}^d$, $d\geq 2$,  $-\infty\leq S<T\leq \infty$, $u: (S,T)\times\Omega \rightarrow \mathbb{R}^d$ denotes the velocity field, $p:(S,T)\times\Omega \rightarrow \mathbb{R}$ denotes the associated pressure, and $f$ and $g$ are a given vector field and a function defined on $(S,T)\times\Omega$, respectively. The variable viscosity coefficients $a^{ij}$ satisfy the following uniform ellipticity conditions: there exists a constant $\nu \in (0,1)$ such that
\begin{equation}\label{eq:elliptic}
\nu|\xi|^2 \leq a^{ij}(t,x) \xi_i\xi_j,\quad |a^{ij}(t,x)|\leq \nu^{-1},\quad \text{for all } i,j\in \{1,\dots,d\}
\end{equation}
for all $\xi=(\xi_1,\dots,\xi_d)\in \mathbb{R}^d$ and $(t,x) \in \mathbb{R}^{d+1}$. We follow the Einstein summation convention for repeated indices unless it is specified otherwise.

The Stokes equation \eqref{eq:Stokes} with variable viscosity coefficients arises in the study of non-Newtonian fluids that exhibit time-dependent shear-thinning properties (see, e.g., \cite{BGL16}). Additionally, it naturally appears in the context of Stokes equations on manifolds (see, e.g., \cite{DM04, MT01}).

As a natural extension of regularity theory for elliptic and parabolic equations, one can ask the validity of the following  local regularity estimates for \eqref{eq:Stokes} under suitable boundary conditions: for $s,q\in (1,\infty)$, $r\in (0,R)$, and $z_0=(0,x_0)$, where $x_0 \in \partial\Omega$, do we have
\begin{equation}\label{eq:Hessian-curved}
\norm{D^2 u}{\Leb{s,q}(\widehat{Q}_r(z_0))}\leq N\left(\norm{u}{\Leb{s,q}(\widehat{Q}_R(z_0))}+\norm{f}{\Leb{s,q}(\widehat{Q}_R(z_0))}+\norm{Dg}{\Leb{s,q}(\widehat{Q}_R(z_0))} \right),
\end{equation}
where $\widehat{Q}_r(z_0)=(\mathbb{R}\times\Omega) \cap Q_r(z_0)$ and $N$ is a constant?

While there has been significant recent progress in obtaining interior estimates for \eqref{eq:Stokes} (see, e.g., \cite{CSTY09, J13, W15, HLW14, DP21, DL22, DK23}), deriving boundary estimates of the type \eqref{eq:Hessian-curved} remains challenging due to the non-local effects of pressure near the boundary. On the one hand, when $\Omega = \mathbb{R}^d_+$, Chang and Kang \cite{CK20} recently demonstrated that such boundary estimates are not feasible under no-slip boundary conditions for \eqref{eq:Stokes}, even when $g = 0$. On the other hand, if additional regularity is assumed for the pressure $p$, the following local regularity estimate can be obtained under no-slip boundary conditions:  
\begin{equation}\label{eq:local-regularity-estimate}
\begin{aligned}
&\norm{\partial_t u}{\Leb{s,q}(\widehat{Q}_r(z_0))}+\norm{D^2 u}{\Leb{s,q}(\widehat{Q}_r(z_0))}+\norm{\nabla p}{\Leb{s,q}(\widehat{Q}_r(z_0))}\\
&\leq N\left(\norm{u}{\Leb{s,q}(\widehat{Q}_R(z_0))}+\norm{p}{\Leb{s,q}(\widehat{Q}_R(z_0))}+\norm{f}{\Leb{s,q}(\widehat{Q}_R(z_0))}+\norm{g}{\Sob{0,1}{s,q}(\widehat{Q}_R(z_0))}\right)
\end{aligned}
\end{equation}
for some constant \( N > 0 \). When \( a^{ij} = \delta^{ij} \) and \( g = 0 \), such results have been studied by Seregin \cite{S00} near flat portions of the boundary and by Vyalov and Shilkin \cite{VS13} near curved boundary segments. For further details, see the exposition by Seregin and Shilkin in \cite[Theorem 2.2]{SS14}.

Another natural boundary condition in fluid mechanics is the slip boundary condition, which consists of impenetrability condition and the slip condition:
\begin{equation}\label{eq:Navier-curved}
u\cdot n=0\quad\text{and}\quad [2\mathbb{D}(u)n+\mathcal{A}u]_{\mathrm{tan}}=0\quad \text{on } \partial\Omega,
\end{equation}
where $n$ is the unit normal vector on $\partial\Omega$ and $[f]_{\mathrm{tan}}=f-(f\cdot n)n$. Here $\mathbb{D}(u)$ is the deformation tensor defined by $(2\mathbb{D}(u))_{ij}=D_i u^j + D_j u^i$ and $\mathcal{A}=\mathcal{A}(t,x)$ is a $2$-tensor on the boundary that describes the friction near the boundary. In the local coordinate of the boundary, $\mathcal{A}$ can be written in matrix form. Since $u\cdot n=0$, $\mathcal{A}u$ lies in the tangent space. 

When $\mathcal{A}=\alpha I_d$, $\alpha\geq 0$, the condition \eqref{eq:Navier-curved} is the classical Navier boundary condition introduced by Navier in 1827. Slip boundary conditions have recently garnered significant attention in computational and experimental studies, particularly in the context of microfluidics (see, e.g., \cite{BB99, LS03, TT97}). For a comprehensive survey on experimental studies concerning slip effects in fluids, we refer the reader to \cite{NEBJ05}.

Navier boundary conditions have also been extensively studied by researchers investigating vanishing viscosity limit problems for incompressible Navier-Stokes equations in bounded domains. Notably, Kelliher \cite{K06}, Xiao-Xin \cite{XX07}, and Gie-Kelliher \cite{GK12} validated the vanishing viscosity limit under the following boundary condition:  
\begin{equation*}\label{eq:Xiao-Xin-vorticity}
 u \cdot n =0\quad\text{and}\quad \omega \times n=0\quad \text{on } \partial\Omega,\quad \Omega \subset \mathbb{R}^3,
\end{equation*}  
where \(\omega\) denotes the vorticity of \(u\).  This boundary condition can be extended to general dimensions as follows:  
\begin{equation}\label{eq:general-Lions-boundary}
u\cdot n=0\quad\text{and}\quad \bm{\omega}n=0\quad\text{on } \partial\Omega,
\end{equation}  
where \(\bm{\omega} = \frac{1}{2}((\nabla u)-(\nabla u)^T)\) is the skew-symmetric part of the velocity gradient \(\nabla u\), often referred to as the vorticity matrix.  
In fact, when $\mathcal{A}$ is the shape operator of $\partial\Omega$, then it can be shown that \eqref{eq:Navier-curved} and \eqref{eq:general-Lions-boundary} are equivalent (see Appendix \ref{sec:principal}).


When $\Omega=\mathbb{R}^d_+$ and $\mathcal{A}=\alpha I_d$, \eqref{eq:Navier-curved} is reduced to
\begin{equation}\label{eq:Navier-flat}
   u^d=0\quad\text{and}\quad -D_d u^k+\alpha u^k=0,\quad k=1,\dots,d-1.
\end{equation}
In particular, when $\alpha=0$, the boundary condition is called the Lions boundary conditions. 
In this case, Dong-Kim-Phan \cite{DKP22} proved that estimate \eqref{eq:Hessian-curved} holds if $(u,p)\in \tilde{W}^{1,2}_{s,q}(\widehat{Q}_R(z_0))^d\times \Sob{0,1}{1}(\widehat{Q}_R(z_0))$ is a strong solution of \eqref{eq:Stokes} in $\widehat{Q}_R(z_0)$ with the Lions boundary conditions and $a^{ij}$ has small mean oscillation with respect to $x$ in small cylinders. Here $\tilde{W}^{1,2}_{s,q}(U) = \{ u : u, Du, D^2 u \in \Leb{s,q}(U), \partial_t u \in L_1(U)\}$. Later, the authors \cite{DK23} relaxed this a priori assumption to $u\in \tilde{W}^{1,2}_{q_0}(\widehat{Q}_R(z_0))^d$ for some $q_0 \in (1,\infty)$. For the case $\alpha>0$, Chen-Liang-Tsai \cite{CLT23} obtained gradient estimates for nonstationary Stokes equations when $a^{ij}=\delta^{ij}$, $f=\Div \mathbb{F}$, and $g=0$. By assuming suitable regularity on the pressure, they also obtained local regularity estimates \eqref{eq:local-regularity-estimate} for weak solutions. Under the same assumption, very recently, the same authors \cite{CLT24} constructed a very weak solution $u$ to nonstationary Stokes equations satisfying \eqref{eq:Navier-flat}, $f=g=0$, and 
\[ \norm{Du}{\Leb{\infty}(Q_1^+)}+\norm{u}{\Leb{\infty}(Q_1^+)}<\infty, \quad \text{but } \norm{D^2 u}{\Leb{q}(Q_{1/2}^+)}=\infty.\]
when $\alpha>0$. In other words, we need to impose additional conditions in order to estimate the Hessian of solutions.

All of the aforementioned results \cite{DP21,DL22,DKP22,DK23,CLT23} were obtained near a flat portion of the boundary. Thus, it is natural to ask whether the estimates \eqref{eq:Hessian-curved} and \eqref{eq:local-regularity-estimate} remain valid near a curved portion of the boundary.

The purpose of this paper is twofold. First, we derive local boundary mixed-norm estimates for the Hessian of solutions in a neighborhood of a curved boundary when \( u \in \Sob{1,2}{q_0} \) satisfies the boundary condition \eqref{eq:general-Lions-boundary} for some \( q_0 \in (1,\infty) \). Second, when $\mathcal{A}$ is sufficiently smooth, we establish a priori local regularity estimates for solutions under the generalized Navier boundary conditions \eqref{eq:Navier-curved}.

Our results are the first to provide boundary \(\Leb{s,q}\)-estimates for Stokes equations near a curved portion of the domain boundary under slip boundary conditions. The first main theorem (Theorem \ref{thm:A}) reveals that, to obtain boundary Hessian \(\Leb{s,q}\)-estimates, there must be a close relationship between the slip tensor \(\mathcal{A}\) and the geometry of the boundary. This extends previous results \cite{CLT23,DK23}, which only considered flat portions of the boundary. Furthermore, we emphasize that the second main theorem (Theorem \ref{thm:B}) does not contradict the counterexample in \cite{CLT24}, as we impose an additional regularity assumption on the pressure and the temporal regularity of the velocity field. Notably, our results are novel even when \( a^{ij} = \delta^{ij} \). Precise statements of the results will be presented in the following section.

\subsubsection*{Outline of the proofs}
Let us outline the proofs of the main theorems. To prove Theorem \ref{thm:A}, a natural approach is to use a standard flattening map and apply a regularity result for the Stokes system with lower-order terms in the half-space to the transformed equation. However, in general, if \( u \) satisfies \eqref{eq:Stokes} and \eqref{eq:Navier-curved}, the transformed equation introduces the derivative of \( p \) on the right-hand side (see \eqref{eq:usual-flattening}). Additionally, the boundary condition \eqref{eq:Navier-curved} becomes significantly more complicated (see \eqref{eq:ud-expression} and \eqref{eq:ud-k}). As a result, it becomes impossible to estimate the Hessian of solutions without imposing additional regularity assumptions on the pressure \( p \) in the corresponding functional space.

To address this issue, we carefully select a flattening map that preserves the impermeability condition (see \eqref{eq:transformation-1}). Under the boundary condition \eqref{eq:general-Lions-boundary}, this allows us to reduce the problem to analyzing the following perturbed Stokes equations with Lions boundary conditions:

\begin{equation}\label{eq:perturbed-Stokes}
\left\{
\begin{aligned}
\partial_t {u} - {a}^{ij}(t,x)D_{ij} {u} - {b}^i(t,x)D_i {u}-c(t,x){u}+\nabla {p} &={f},\\
\Div(\widehat{A}{u})&={g},
\end{aligned}
\right.
\end{equation}
and 
\begin{equation}\label{eq:Lions-boundary}
{u}^d=0\quad \text{and}\quad D_d {u}^k=0,\quad k=1,\dots,d-1,
\end{equation}
where \(\widehat{A}\) is sufficiently close to the identity (see Assumption \ref{assump:delta-N}), and \(b\) and \(c\) are bounded. 
We then follow a strategy developed by the authors in \cite{DK23}, adapted to this setting. First, we establish a priori local boundary Hessian estimates for \eqref{eq:perturbed-Stokes} in a small half-cylinder \(Q_R^+\) (see Theorem \ref{thm:regularity-a-priori-local}). Specifically, if \((u,p) \in \Sob{1,2}{s,q}(Q_R^+)^d \times \Sob{0,1}{1}(Q_R^+)\) is a strong solution of \eqref{eq:perturbed-Stokes} in \(Q_R^+\) with Lions boundary conditions \eqref{eq:Lions-boundary}, for some \(f \in \Leb{s,q}(Q_R^+)^d\) and \(g \in \Sob{0,1}{s,q}(Q_R^+)\), we have the following a priori estimate:  
\[
\norm{D^2 u}{\Leb{s,q}(Q_r^+)}\leq N \left( \norm{u}{\Leb{s,1}(Q_R^+)} + \norm{f}{\Leb{s,q}(Q_R^+)} + \norm{Dg}{\Leb{s,q}(Q_R^+)} \right)
\]
for some constant \(N > 0\). 

To relax the assumption \(u \in \Sob{1,2}{s,q}(Q_R^+)^d\), we utilize the boundary condition \eqref{eq:Lions-boundary} to extend \((u,p)\) to the entire cylinder \(Q_R\). For a small \(\varepsilon > 0\), we then perform a space-time mollification of \eqref{eq:perturbed-Stokes}, yielding the following equations: 
\[
\partial_t u^{(\varepsilon)} - \mathcal{L}u^{(\varepsilon)} + \nabla p^{(\varepsilon)} = f^{(\varepsilon)} + {f}^\varepsilon_e \quad \text{in } Q_R,
\]
and
\[
\Div(\widehat{A}u^{(\varepsilon)}) = g^{(\varepsilon)} + g^\varepsilon_e \quad \text{in } Q_R,
\]
where
\begin{align*}
\mathcal{L}u &= a^{ij}(t,x)D_{ij} u + b^i(t,x) D_i u + c(t,x)u,\\
{f}^\varepsilon_e &= \left[(\mathcal{L}u)^{(\varepsilon)} - \mathcal{L}u^{(\varepsilon)}\right]1_{Q_R},\quad
{g}^\varepsilon_e = \Div\left[(\widehat{A}u^{(\varepsilon)} - (\widehat{A}u)^{(\varepsilon)})\psi\right],
\end{align*}
and \(a^{ij}\), \(b^i\), and \(c\) are appropriately extended to the whole space.

We decompose $u^{(\varepsilon)}=u_1^\varepsilon+u_2^\varepsilon$, where $u_1^\varepsilon$ is obtained by using a $\Sob{1,2}{q_0}$ solvability result (Theorem \ref{thm:solvability-modified-divergence}) for \eqref{eq:perturbed-Stokes} with $f_e^\varepsilon$ and $g^\varepsilon_e$ instead of $f$ and $g$, where $f_e^\varepsilon,g^\varepsilon_e\rightarrow 0$ in an appropriate sense. Another consequence of the weighted solvability result is that we can improve the regularity of $u_2^\varepsilon=u^{(\varepsilon)}-u_1^\varepsilon$ (Lemma \ref{lem:regularity-lemma-Stokes}), allowing us to apply a priori local boundary estimate to $u_2^\varepsilon$. Then the desired estimates for the perturbed system follow by a compactness argument. Finally, a standard covering argument provides the desired regularity estimates for the original system.

To prove Theorem \ref{thm:B}, we use a standard flattening map to transform the original problem into another perturbed Stokes equation \eqref{eq:transformed-A} satisfying the boundary condition \eqref{eq:transformed-B}. Unlike the previous case, the transformed solution does not satisfy the Lions boundary condition \eqref{eq:Lions-boundary} due to the absence of the special structure on $\mathcal{A}$. To address this, we first correct these error terms by using a solution $v$ of the boundary value problems of the heat equation on the half-space so that their difference satisfies the Lions boundary conditions. The next natural step is to apply the maximal $\Leb{s,q}$-regularity results to get a local estimate. However, this correction introduces an additional error term involving $\partial_t u$ on the right-hand side of the inequality. To eliminate this term, we employ a classical idea by Agmon (Lemma \ref{lem:lambda-estimate}; see also \cite[Lemma 5.5]{K07}) to absorb the norm of $\partial_t u$ via an iteration argument. Then the desired result follows by a covering argument.

\subsubsection*{Organizations}
This paper is organized into six sections and an appendix. The notation and main results are introduced in the next section. In Section \ref{sec:prelim}, we define the function spaces and trace operators, and collect key results on the solvability of boundary value problems for the heat equation and Stokes equations with variable coefficients in the half-space, as well as relevant interpolation inequalities. Additionally, we introduce a boundary-flattening map that preserves the impermeability condition, which is crucial for proving Theorem \ref{thm:A}. In Section \ref{sec:perturbed-Stokes}, we establish solvability and local boundary regularity estimates for perturbed Stokes equations. The proofs of Theorems \ref{thm:A} and \ref{thm:B} are presented in Sections \ref{sec:5} and \ref{sec:6}, respectively.

In Appendix \ref{sec:principal}, we connect the Navier boundary conditions \eqref{eq:Navier-curved} with \eqref{eq:general-Lions-boundary} when $\mathcal{A}$ becomes the shape operator of $\partial\Omega$.

\section{Notation and main results}

\subsection{Notation}
By $N=N(p_1,\dots,p_k)$, we denote a generic positive constant depending only on the parameters $p_1,\dots,p_k$. 
For $x_0=(x_0',x_{0d}) \in \mathbb{R}^{d-1}\times \mathbb{R}$, we write 
\[  B_\rho(x_0)=\{ y \in \mathbb{R}^d : |x_0-y|<\rho\}\]
and $z_0=(t_0,x_0) \in \mathbb{R}^{d+1}$, we write
\[  Q_\rho(z_0) = (t_0-\rho^2,t_0)\times B_\rho(x_0)\]
the parabolic cylinder centered at $z_0$ with radius $\rho>0$. We write
\[  B_\rho^+(x_0)=B_\rho(x_0)\cap \mathbb{R}^d_+,\quad Q_\rho^+(z_0)=Q_\rho(z_0)\cap (\mathbb{R}\times\mathbb{R}^d_+),\]
and 
\[  B'_\rho(x_0')=\{ y' \in \mathbb{R}^{d-1} : |x'-y'|<\rho\}.\]
For $z_0=(t_0,x_0) \in \mathbb{R}\times \partial\Omega$, we write $\widehat{B}_\rho(x_0)=B_\rho(x_0)\cap \Omega$ and  $\widehat{Q}_\rho(z_0)=Q_\rho(z_0)\cap (\mathbb{R}\times \Omega)$. 
When $x_0=0$ and $z_0=(0,0)$, we write $B_\rho=B_\rho(0)$ and $Q_\rho=Q_\rho(0,0)$. Similarly, we write $B_\rho'= B_\rho'(0)$, $B_\rho^+=B_\rho^+(0)$, $Q_\rho^+=Q_\rho^+(0,0)$, $\widehat{B}_\rho =\widehat{B}_\rho(0)$, and $\widehat{Q}_\rho = \widehat{Q}_\rho(0,0)$. For $U \subset \mathbb{R}^{d+1}$, let $C_0^\infty (U)$ denote $C^\infty$-functions defined on $U$ that have compact supports in $U$. 

For $U=\Gamma \times \Omega \subset \mathbb{R}\times\mathbb{R}^d$ and $1\leq s,q<\infty$, we write 
\[  \norm{h}{\Leb{s,q}(U)}:=\left[\int_\Gamma \left(\int_\Omega |h(t,x)|^q \myd{x} \right)^{s/q} dt \right]^{1/s} \]
and we denote the mixed-norm Lebesgue space 
\[  \Leb{s,q}(U):=\{ h : \norm{h}{\Leb{s,q}(U)}<\infty \}.\]
Similarly, one can define $\Leb{s,q}(U)$ either $s=\infty$ or $q=\infty$. The parabolic Sobolev spaces of mixed type are defined by
\[   \Sob{1,2}{s,q}(U) := \{ u : u, Du, D^2 u , \partial_t u \in \Leb{s,q}(U)\}.\]
We also define
\[  \Sob{0,k}{s,q}(U) := \{ u : u, Du, ..., D^k u \in \Leb{s,q}(U)\},\quad k=1,2.\]
When $s=q$, we simply write $\Leb{q}(U)=\Leb{q,q}(U)$, $\Sob{1,2}{q}(U)=\Sob{1,2}{q,q}(U)$, and $\Sob{0,k}{q}(U)=\Sob{0,k}{q,q}(U)$. 

 For $A\subset \mathbb{R}^d$, we write the spatial average of $f$ over $A$ by 
\[  [f]_A(t):=\frac{1}{|A|}\int_A f(t,x)\myd{x},\]
where $|A|$ denotes the $d$-dimensional Lebesgue measure of $A$. Similarly, we write the total average of $f$ over $\mathcal{D}\subset \mathbb{R}\times\mathbb{R}^d$ by 
\[  (f)_{\mathcal{D}}=\fint_{\mathcal{D}} f \myd{x}dt:=\frac{1}{|\mathcal{D}|}\int_{\mathcal{D}} f(t,x)\myd{x}dt,\]
where $|\mathcal{D}|$ denotes the $(d+1)$-dimensional Lebesgue measure of $\mathcal{D}$.

\subsection{Main results}
We impose the following assumption on the viscosity coefficients.
\begin{assumption}[$\gamma$]\label{assump:VMO}
There exists $R_0>0$ such that 
\[ \fint_{Q_r(t_0,x_0)} |a^{ij}(t,x)-[a^{ij}]_{B_r(x_0)}(t)|\myd{x}dt\leq \gamma \]
for any $(t_0,x_0)\in\mathbb{R}^{d+1}$, $0<r<R_0$, and for all $i,j=1,\dots,d$.
\end{assumption}

To state our assumption on the boundary of the domain, we say that $\phi \in C^{2,1}(\mathbb{R}^d)$ if $\phi \in C^2(\mathbb{R}^d)$ and $D_{ij} \phi$ is Lipschitz continuous on $\mathbb{R}^d$ for all $i,j$. The norm is defined by 
\[ \norm{\phi}{C^{2,1}(\mathbb{R}^d)}=\norm{\phi}{\Leb{\infty}(\mathbb{R}^d)}+\norm{D\phi}{\Leb{\infty}(\mathbb{R}^d)}+\norm{D^2\phi}{\Leb{\infty}(\mathbb{R}^d)}+[D^2\phi]_{C^{0,1}(\mathbb{R}^d)},\]
where $[D^2\phi]_{C^{0,1}(\mathbb{R}^d)}$ denotes the Lipschitz norm of $D^2 \phi$ on $\mathbb{R}^d$.

\begin{assumption}[$\theta,R_0^*$]\label{assump:domain}
There exists $M>0$ such that for any $x_0 \in \partial\Omega$, there exist a $C^{2,1}$-function $\phi:\mathbb{R}^{d-1}\rightarrow\mathbb{R}$ and the Cartesian coordinate system so that $\phi(0)=0$, $\nabla \phi(0)=0$, 
\[  \Omega \cap B_{R_0^*}(x_0)=\{ x \in B_{R_0^*}(x_0) : x_d >\phi(x')\}, \] 
and 
\[   \norm{\nabla \phi}{\Leb{\infty}(B_{R_0^*}'(x_0'))} \leq \theta \quad \text{and}\quad \norm{\phi}{C^{2,1}(\mathbb{R}^{d-1})}\leq M\]
in a coordinate system.
\end{assumption}
\begin{remark}
Since $\nabla \phi$ is continuous and $\nabla \phi(x_0)=0$, for any $\theta>0$, there exists $R_0^*>0$ such that $|\nabla \phi(x')|<\theta$ for all $(x',x_d)\in B_{R_0^*}(x_0)$.
\end{remark}

The first result concerns the local boundary Hessian estimates of solutions to Stokes equations under the boundary condition \eqref{eq:general-Lions-boundary}. 

\begin{theorem}\label{thm:A}
Let $s,q_0,q\in (1,\infty)$, $r\in (0,R)$, $x_0\in\partial \Omega,$ and $z_0=(0,x_0)$. Then there exists a constant $\gamma>0$ depending on $d$, $s$, $q$, $q_0$, $\nu$ such that under Assumption \ref{assump:VMO} $(\gamma)$, if $(u,p) \in {W}^{1,2}_{q_0}(\widehat{Q}_R(z_0))^d \times \Sob{0,1}{1}(\widehat{Q}_R(z_0))$ is a strong solution of \eqref{eq:Stokes} in $\widehat{Q}_R(z_0)$ with the boundary condition \eqref{eq:general-Lions-boundary} on $Q_R(z_0)\cap (\mathbb{R}\times \partial\Omega)$ for some $f\in \Leb{s,q}(\widehat{Q}_R(z_0))^d$ and $g\in \Sob{0,1}{s,q}(\widehat{Q}_R(z_0))$, then $D^2 u \in \Leb{s,q}(\widehat{Q}_r(z_0))$. Moreover, we have
\[   \norm{D^2 u}{\Leb{s,q}(\widehat{Q}_r(z_0))}\leq N\left( \norm{u}{\Leb{s,1}(\widehat{Q}_R(z_0))} +  \norm{f}{\Leb{s,q}(\widehat{Q}_R(z_0))}+ \norm{g}{\Sob{0,1}{s,q}(\widehat{Q}_R(z_0))}\right)\]
for some constant $N=N(d,s,q,\nu,\Omega,r,R,R_0)>0$.
\end{theorem}

The second result concerns a priori local boundary estimates for Stokes equations for general slip tensor $\mathcal{A} \in C^{1,2}_{par}$, i.e., $\partial_t\mathcal{A}$, $\mathcal{A}$, $D_x\mathcal{A}$, $D^2_x\mathcal{A}$ are bounded and continuous. 

\begin{theorem}\label{thm:B}
Let $\mathcal{A} \in C^{1,2}_{par}(\mathbb{R}\times\partial\Omega)$, $s,q \in (1,\infty)$, $r\in (0,R)$, $x_0\in\partial \Omega,$ and $z_0=(0,x_0)$. Then there exists a constant $\gamma>0$ depending on $d$, $s$, $q$, $\nu$ such that under Assumption \ref{assump:VMO} $(\gamma)$, if $(u,p) \in {W}^{1,2}_{s,q}(\widehat{Q}_R(z_0))^d \times \Sob{0,1}{s,q}(\widehat{Q}_R(z_0))$ is a strong solution of \eqref{eq:Stokes} in $\widehat{Q}_R(z_0)$ with the boundary condition \eqref{eq:Navier-curved} on $Q_R(z_0)\cap (\mathbb{R}\times \partial\Omega)$ for some $f\in \Leb{s,q}(\widehat{Q}_R(z_0))^d$ and $g\in \Sob{0,1}{s,q}(\widehat{Q}_R(z_0))$ with $\partial_t g =\Div G$ in $\widehat{Q}_R(z_0)$ for some $G\in \Leb{s,q}(\widehat{Q}_R(z_0))^d$, then we have
\begin{align*}
&\norm{\partial_t u}{\Leb{s,q}(\widehat{Q}_r(z_0))}+\norm{D^2{u}}{\Leb{s,q}(\widehat{Q}_r(z_0))}+\norm{\nabla {p}}{\Leb{s,q}(\widehat{Q}_r(z_0))}\\
&\leq N\left(\norm{{u}}{\Leb{s,q}(\widehat{Q}_R(z_0))}+\norm{{p}}{\Leb{s,q}(\widehat{Q}_R(z_0))}+\norm{{f}}{\Leb{s,q}(\widehat{Q}_R(z_0))}\right.\\
&\relphantom{=}\left.\qquad+\norm{{g}}{\Sob{0,1}{s,q}(\widehat{Q}_R(z_0))}+\norm{G}{\Leb{s,q}(\widehat{Q}_R(z_0))}\right)
\end{align*}
for some constant $N=N(d,s,q,\nu,\mathcal{A},\Omega,r,R,R_0)>0$.
\end{theorem}

\begin{remark}\leavevmode
\begin{enumerate}
\item When $x_0$ lies on a flat portion of the boundary, Theorem \ref{thm:A} was obtained by the authors \cite[Theorem 2.10]{DK23}. Unlike the flat case, Theorem \ref{thm:A} requires us to assume that $\partial_t u \in \Leb{q_0}(\hat{Q}_R(z_0))$ for some $q_0>1$. See Remark \ref{rem:u-t-estimate}.
\item Theorem \ref{thm:B} concerns a priori local boundary mixed norm estimate. When $a^{ij}=\delta^{ij}$, $g=0$, and $\mathcal{A}=\alpha I_d$, $\alpha>0$ a constant, Chen-Liang-Tsai \cite[Theorem 1.2]{CLT23} proved the estimate on a neighborhood of a flat portion of the boundary. Moreover, they proved that a weak solution becomes strong. It would be interesting to see if one could prove spatial smoothing of solutions when $f \in \Leb{s,q}(\widehat{Q}_R(z_0))^d$, $g\in \Sob{0,1}{s,q}(\widehat{Q}_R(z_0))$ with $\partial_t g =\Div G$ in $\widehat{Q}_R(z_0)$ for some $G\in \Leb{s,q}(\widehat{Q}_R(z_0))^d$. For the no-slip boundary case, such an estimate was obtained by Vyalov-Shilkin \cite{VS13} when $a^{ij}=\delta^{ij}$ and $g=0$.
\end{enumerate} 
\end{remark}

\section{Preliminaries}\label{sec:prelim}
This section consists of five subsections. We first introduce function spaces including weighted Sobolev spaces, anisotropic Bessel potential spaces, and Triebel-Lizorkin spaces. We also recall interpolation inequalities that will be used in this paper. In Section \ref{subsec:trace-operator}, we state the spatial trace theorem and the solvability of the heat equation with either Dirichlet or Neumann boundary conditions. In Section \ref{subsec:boundary-flattening}, we introduce the boundary flattening map which preserves the impermeability condition. In Section \ref{subsec:Stokes-variable}, we recall the solvability of Stokes equations with variable coefficients on the half-space under the Lions boundary conditions. Finally, in the last subsection, we recall local interior and boundary Hessian estimates for Stokes equations with variable coefficients.

\subsection{Function spaces}

A function $w$ is a \emph{weight} on $\mathbb{R}^d$ if $w$ is nonnegative and $w>0$ a.e. on $\mathbb{R}^d$. For $q\in (1,\infty)$, we write $w\in A_q(\mathbb{R}^d,dx)$ if
\[ [w]_{A_q} :=\sup_{x_0\in \mathbb{R}^d, r>0} \left(\fint_{B_r(x_0)} w \myd{x} \right)\left(\fint_{B_r(x_0)} w^{-1/(q-1)}\myd{x}\right)^{q-1} <\infty.
\]

For $s,q \in (1,\infty)$, $K\geq 1$, and a function $w$ on $\mathbb{R}^{d+1}$, we write $[w]_{A_{s,q}} \leq K$ if there exist weights $w_1$ on $\mathbb{R}^d$ and $w_2$ on $\mathbb{R}$ such that
\begin{equation*}
  w(t,x)=w_1(x)w_2(t)\quad \text{and}\quad   [w_1]_{A_q}, [w_2]_{A_s}\leq K.
\end{equation*}

For $s,q\in (1,\infty)$, $-\infty\leq S<\infty$, $-\infty<T\leq \infty$, and the weight $w \in A_{s,q}$, we define
\begin{align*}
\norm{f}{\Leb{s,q,w}((S,T)\times\Omega)}&:=\left(\int_S^T\left(\int_\Omega |f|^q w_1 \myd{x}\right)^{s/q} w_2 \myd{t}\right)^{1/s}
\end{align*}
and
\[
  \Leb{s,q,w}((S,T)\times\Omega):=\{ f : \norm{f}{\Leb{s,q,w}((S,T)\times\Omega)}<\infty \}.
\]
Similarly, for  $s,q\in(1,\infty)$ and $w\in A_{s,q}$, we define weighted parabolic Sobolev spaces
\begin{align*}
{\Sob{0,1}{s,q,w}((S,T)\times\Omega)}&:=\{ u : u, Du \in \Leb{s,q,w}((S,T)\times \Omega)\},\\
    \Sob{1,2}{s,q,w}((S,T)\times\Omega)&:=\{ u : u, Du, D^2 u, \partial_t u \in\Leb{s,q,w}((S,T)\times\Omega)\}
\end{align*}
with the norm
\begin{align*}
\norm{u}{\Sob{0,1}{s,q,w}((S,T)\times\Omega)}&:=\norm{u}{\Leb{s,q,w}((S,T)\times\Omega)}+\norm{Du}{\Leb{s,q,w}((S,T)\times\Omega)},\\
   \norm{u}{\Sob{1,2}{s,q,w}((S,T)\times\Omega)}&:=\norm{\partial_t u}{\Leb{s,q,w}((S, T)\times\Omega)}+\sum_{k=0}^2\norm{D^k u}{\Leb{s,q,w}((S,T)\times\Omega)}.
\end{align*}

For $b\in \mathbb{R}$ and $s,q\in (1,\infty)$, let $\fSob{b/2,b}{s,q}(\mathbb{R}\times \mathbb{R}^d)$ and $F^{b/2,b}_{s,q}(\mathbb{R}\times \mathbb{R}^d)$ be the anisotropic Bessel potential spaces and Triebel-Lizorkin spaces which can be defined via Littlewood-Paley projections (see e.g. \cite{ST87}). We also write  $[f]_{F^{b/2,b}_{s,q}(\mathbb{R}\times\mathbb{R}^d)}$ the seminorm of the anisotropic Triebel-Lizorkin spaces.

When $b> 0$, the anisotropic Bessel potential spaces $\fSob{b/2,b}{s,q}(\mathbb{R}\times \mathbb{R}^d)$ can be characterized as 
\begin{equation}\label{eq:anisotropic-Bessel}
 \fSob{b/2,b}{s,q}(\mathbb{R}\times \mathbb{R}^d)=\Leb{s}(\mathbb{R};\fSob{b}{q}(\mathbb{R}^d)) \cap \fSob{b/2}{s}(\mathbb{R};\Leb{q}(\mathbb{R}^d)).
\end{equation}
Here $\fSob{b}{q}$ is the isotropic Bessel potential spaces. For a domain $\Omega$ in $\mathbb{R}^d$, we denote by $\fSob{b/2,b}{s,q}(\mathbb{R}\times\Omega)$ the Bessel potential spaces on $\mathbb{R}\times\Omega$, where 
\[
  \fSob{b/2,b}{s,q}(\mathbb{R}\times\Omega)= \{ Rf : f \in \fSob{b/2,b}{s,q}(\mathbb{R}\times \mathbb{R}^d) \},
\]
where $Rf$ denotes the restriction of $f$ on $\mathbb{R}\times \Omega$. The norm is defined by 	
\[ \norm{f}{\fSob{b/2,b}{s,q}(\mathbb{R}\times\Omega)} =\inf\{ \norm{g}{\fSob{b/2,b}{s,q}(\mathbb{R}\times\mathbb{R}^d)} : f=Rg,\, g\in \fSob{b/2,b}{s,q}(\mathbb{R}\times \mathbb{R}^d)\}.\]
 Also, it can be characterized as 
\[ \fSob{b/2,b}{s,q}(\mathbb{R}\times\Omega) = \Leb{s}(\mathbb{R};\fSob{b}{q}(\Omega)) \cap \fSob{b/2}{s}(\mathbb{R};\Leb{q}(\Omega)), \]
where $\fSob{b}{q}(\Omega)$ is defined in a similar manner (see \cite[Appendix]{G18}).

The following mapping property can be found in \cite[Appendix]{G18}.
\begin{lemma}\label{lem:derivative-map}
Let $s,q\in (1,\infty)$ and $2\kappa \in \{0,1,2\}$. Then $\partial_{x_j}$, $j=1,\dots,d$, are continuous from $\fSob{\kappa,2\kappa}{s,q}(\mathbb{R}\times \mathbb{R}^d_+)$ to $\fSob{\kappa-1/2,2\kappa-1}{s,q}(\mathbb{R}\times \mathbb{R}^d_+)$.
\end{lemma}

We recall the following interpolation inequalities. 

\begin{lemma}\label{lem:interpolation-inequality}
Let $1\leq p\leq q<\infty$, $q>1$, and $R>0$. Suppose that $w\in A_q$ with $[w]_{A_q}\leq K_0$ and $\varepsilon>0$. Then there exist constants $N_1=N_1(d,p,q)>0$ and $N_2=N_2(d,q,K_0)>0$ such that 
\begin{align}
\norm{Du}{\Leb{q}(B_R^+)}&\leq N_1 \left(\norm{D^2 u}{\Leb{q}(B_R^+)}^{\theta} \norm{u}{\Leb{p}(B_R^+)}^{1-\theta} +  R^{-1+d(1/q-1/p)} \norm{u}{\Leb{p}(B_R^+)}\right),\label{eq:domain-interpolation}\\
  \norm{Du}{\Leb{q,w}(\mathbb{R}^d)}&\leq \varepsilon\norm{D^2 u}{\Leb{q,w}(\mathbb{R}^d)}+\frac{N_2}{\varepsilon} \norm{u}{\Leb{q,w}(\mathbb{R}^d)},\label{eq:weighted-interpolation}
\end{align}
where $\theta = (1/p+1/d-1/q)/(1/p+2/d-1/q)$ and $\theta \in [1/2,1)$.
\end{lemma}
\begin{proof}
The inequality \eqref{eq:domain-interpolation} is well known. The inequality \eqref{eq:weighted-interpolation} can be proved by using \cite[Lemma 3.8]{DKr19}  and Young's inequality.
\end{proof}

\begin{remark}
If $p=1$, then by Young's inequality, for every $\varepsilon>0$, we have 
\begin{equation}\label{eq:Du-interpolation-L1}
 \norm{Du}{\Leb{q}(B_R^+)}\leq \varepsilon \norm{D^2u}{\Leb{q}(B_R^+)} +N\left(R^{-1+d(1/q-1)} + \varepsilon^{-\theta/(1-\theta)}\right)\norm{u}{\Leb{1}(B_R^+)}
\end{equation}
for some constant $N=N(d,q)>0$.
\end{remark}
 
\begin{lemma}\label{lem:interpolation-half-derivative}
Let $\Omega \in \{\mathbb{R}^d,\mathbb{R}^d_+\}$, $s,q\in (1,\infty)$, and let $\varepsilon>0$. There exists a constant $N=N(d,s,q)>0$ such that 
\[   \norm{u}{\fSob{1/2,1}{s,q}(\mathbb{R}\times\Omega)}\leq \varepsilon \left(\norm{\partial_t u}{\Leb{s,q}(\mathbb{R}\times \Omega)} + \norm{D^2 u}{\Leb{s,q}(\mathbb{R}\times\Omega)}\right) + N\left(1+\frac{1}{\varepsilon}\right)\norm{u}{\Leb{s,q}(\mathbb{R}\times\Omega)} \]
for all $u\in \Sob{1,2}{s,q}(\mathbb{R}\times\Omega)$.
\end{lemma}
\begin{proof}
Recall from \eqref{eq:anisotropic-Bessel} that 
\begin{align*}
\norm{u}{\fSob{1/2,1}{s,q}(\mathbb{R}\times\Omega)} =\norm{u}{\Leb{s,q}(\mathbb{R}\times\Omega)} + \norm{Du}{\Leb{s,q}(\mathbb{R}\times\Omega)} + \norm{u}{\fSob{1/2}{s}(\mathbb{R};\Leb{q}(\Omega))}.
\end{align*}
By the complex interpolation characterization (see e.g. \cite[Proposition 3.7]{MV15}), we have
\[  \fSob{1/2}{s}(\mathbb{R};\Leb{q}(\Omega))=[\Sob{1}{s}(\mathbb{R};\Leb{q}(\Omega)),\Leb{s}(\mathbb{R};\Leb{q}(\Omega))]_{1/2}.\]
Then by Lemma \ref{lem:interpolation-inequality} with $R=\infty$ and Young's inequality, we have
\begin{align*}
\norm{u}{\fSob{1/2,1}{s,q}(\mathbb{R}\times\Omega)}&\leq N \norm{D^2 u}{\Leb{s,q}(\mathbb{R}\times\Omega)}^{1/2} \norm{u}{\Leb{s,q}(\mathbb{R}\times\Omega)}^{1/2}+N\norm{u}{\Sob{1}{s}(\mathbb{R};\Leb{q}(\Omega))}^{1/2}\norm{u}{\Leb{s}(\mathbb{R};\Leb{q}(\Omega))}^{1/2}\\
&\leq \varepsilon \left(\norm{D^2u}{\Leb{s,q}(\mathbb{R}\times\Omega)} + \norm{\partial_t u}{\Leb{s,q}(\mathbb{R}\times\Omega)}\right) + N\left(1+\frac{1}{\varepsilon}\right) \norm{u}{\Leb{s,q}(\mathbb{R}\times\Omega)}
\end{align*}
for some constant $N=N(d,s,q)>0$. This completes the proof of Lemma \ref{lem:interpolation-half-derivative}.
\end{proof}

\subsection{Trace operators and solvability of the inhomogeneous heat equation}\label{subsec:trace-operator}

The following trace theorem can be found in \cite[Lemma 3.5]{DHP07} and \cite[Corollary 3.8]{HL22}.
\begin{theorem}
Let $s,q\in (1,\infty)$ and $2\tau \in \{1,2\}$. Then the trace with respect to $x$
\[ \Tr : \fSob{\tau,2\tau}{s,q}(\mathbb{R}\times \mathbb{R}^d_+)\rightarrow F^{\tau-1/(2q),2\tau-1/q}_{s,q}(\mathbb{R}\times \mathbb{R}^{d-1}) \]
is continuous and surjective.
\end{theorem}

We recall the following solvability results for the heat equation under Dirichlet and Neumann boundary conditions, which can be readily derived from, for instance, \cite{DHP07} and a scaling argument.
\begin{theorem}\label{thm:Dirichlet-heat}
Let $s,q\in(1,\infty)$, $b=2-1/q$, and $\lambda>0$. For any $g\in F^{b/2,b}_{s,q}((0,\infty)\times \mathbb{R}^{d-1})$ and $g(0,\cdot)=0$ when $1/q+2/s<2$, there exists a unique solution $v\in \Sob{1,2}{s,q}((0,\infty)\times \mathbb{R}^d_+)$ satisfying $v(0,\cdot)=0$, 
\[ \partial_t v -\Delta v+\lambda v=0\quad \text{in } (0,\infty)\times \mathbb{R}^d_+,\quad\text{and}\quad v=g\quad \text{on } (0,\infty)\times\mathbb{R}^{d-1}\times\{0\}.\]
Moreover, we  have
\begin{align*}
&\relphantom{=}\norm{\partial_t v}{\Leb{s,q}((0,\infty)\times \mathbb{R}^d_+)}+\sum_{k=0}^2 \lambda^{1-k/2} \norm{D^k v}{\Leb{s,q}((0,\infty)\times \mathbb{R}^d_+)}\\
&\relphantom{=}\leq N \left([g]_{F^{b/2,b}_{s,q}((0,\infty)\times \mathbb{R}^{d-1})}+\lambda^{b/2}\norm{g}{\Leb{s,q}((0,\infty)\times \mathbb{R}^{d-1})} \right)
\end{align*}
for some constant $N=N(d,s,q)>0$. 
\end{theorem}
\begin{theorem}\label{thm:Neumann-heat}
Let $s,q\in (1,\infty)$, $b=1-1/q$, and $\lambda>0$. For any $g\in F^{b/2,b}_{s,q}((0,\infty)\times \mathbb{R}^{d-1})$ and $g(0,\cdot)=0$ when $1/(2q)+1/s<1/2$, there exists a unique solution $v\in \Sob{1,2}{s,q}((0,\infty)\times \mathbb{R}^d_+)$ satisfying $v(0,\cdot)=0$,
\[ \partial_t v -\Delta v+\lambda v=0\quad \text{in } (0,\infty)\times \mathbb{R}^d_+,\quad\text{and}\quad D_d v=g\quad \text{on } (0,\infty)\times\mathbb{R}^{d-1}\times\{0\}.\]
Moreover, we  have
\begin{align*}
&\norm{\partial_t v}{\Leb{s,q}((0,\infty)\times \mathbb{R}^d_+)}+\sum_{k=0}^2 \lambda^{1-k/2} \norm{D^k v}{\Leb{s,q}((0,\infty)\times \mathbb{R}^d_+)}\\
&\leq N \left([g]_{F^{b/2,b}_{s,q}((0,\infty)\times \mathbb{R}^{d-1})}+\lambda^{b/2}\norm{g}{\Leb{s,q}((0,\infty)\times \mathbb{R}^{d-1})} \right)
\end{align*}
for some constant $N=N(d,s,q)>0$. 
\end{theorem}

\subsection{Boundary flattening map preserving impermeability condition}\label{subsec:boundary-flattening}

In this subsection, we introduce a boundary flattening map that plays a crucial role in proving Theorem \ref{thm:A}. A similar flattening map was considered in several papers on the boundary value problem for kinetic equations, e.g., \cite{GHJO20}.

We locally flatten the boundary  $\partial\Omega$ as follows. For $x_0\in \partial\Omega$, there exist $R_0^*>0$, a coordinate system, and a mapping $\phi$ such that $\phi(0)=0$, $\nabla \phi(0)=0$, and 
\[ \Omega_{R_0^*}(x_0):= \Omega \cap B_{R_0^*}(x_0)=\{ x\in B_{R_0^*}(x_0) : x_d>\phi(x')\}.\]

Define a map $\Phi : \mathbb{R}^d \rightarrow \mathbb{R}^d$ by 
\begin{equation}\label{eq:transformation-1}
\left\{\begin{aligned}
x_i&=\Phi^i(y)=y_i -y_d D_i \phi^{(y_d)}(y'),\quad \text{for } i=1,\dots,d-1,\\
x_d&=\Phi^d(y)=y_d + \phi^{(y_d)}(y'),
\end{aligned}
\right.
\end{equation}
where $\phi^{(y_d)}$ is defined as follows: for $y_d \in \mathbb{R}$ and $\zeta \in C_0^\infty(\mathbb{R}^{d-1})$ satisfying $\supp \zeta \subset B_1'$ and $\int_{\mathbb{R}^{d-1}} \zeta \myd{y'}=1$, we define
\[ \phi^{(y_d)}(y')=\int_{\mathbb{R}^{d-1}} \phi(y'-y_d z') \zeta(z')dz'. \]
Note that $\phi^{(0)}(y')=\phi(y')$ and $\phi^{(y_d)}(y')$ is infinitely differentiable in $y$ when $y_d \neq 0$. 

The above regularized boundary function has the following properties:
\begin{lemma}\label{lem:estimates-boundary} Under Assumption \ref{assump:domain} $(\theta,R_0^*)$, there exists a constant $N=N(d,M)>0$ such that 
\begin{enumerate}[label=\textnormal{(\alph*)}]
\item $|D \phi^{(y_d)}(y')|\leq N \theta$ for all $y \in B_{r}$;
\item $|y_d D_{ij} \phi^{(y_d)}(y')|\leq N \theta$ for all $y \in B_{r}$ and for $i,j=1,\dots,d$;
\item $|D^2 \phi^{(y_d)}(y')|+|D^3 \phi^{(y_d)}(y')|\leq N$;
\item $|D^k(y_d D_{id} \phi^{(y_d)}(y'))|\leq N$ for $i=1,\dots,d-1$ and $k=1,2$
\end{enumerate} 
for all $0<r<R^*_0/2$.
\end{lemma}
\begin{proof}
(a) Note that 
\begin{align}
D_i \phi^{(y_d)}(y')&=\int_{\mathbb{R}^{d-1}} (D_i \phi)(y'-y_d z') \zeta(z')dz',\quad i=1,\dots,d-1,\label{eq:Di-phi}\\
D_d \phi^{(y_d)}(y')&=-\int_{\mathbb{R}^{d-1}} (D_{y'} \phi)(y'-y_d z') \cdot (z') \zeta(z')dz'.\label{eq:Dd-phi}
\end{align}
Hence it follows that 
\[  |D \phi^{(y_d)} (y')|\leq N \theta\quad \text{for all } y\in B_{r},\]
where $N=\norm{\zeta}{\Leb{1}}+\norm{|\cdot|\zeta}{\Leb{1}}$. This proves (a).

(b) For $i=1,\dots,d-1$, we have
\begin{equation}\label{eq:yd-identity}   y_d D_{id} \phi^{(y_d)}(y')=\int_{\mathbb{R}^{d-1}} (D_i \phi)(y'-y_d z') [(d-1)\zeta(z')+z' \cdot\nabla \zeta(z')]dz'.  
\end{equation}
Indeed, for $y_d\neq0$, a change of variable and integration by parts give
\begin{align*}
D_{id} \phi^{(y_d)}(y')&=\int_{\mathbb{R}^{d-1}} (D_{il} \phi)(y'-y_d z')(-z_l)\zeta(z')dz'\\
&=\sum_{l=1}^{d-1}\int_{\mathbb{R}^{d-1}} (D_{i}\phi)(y'-z')\frac{\partial}{\partial z_l}\left[-\frac{z_l}{y_d}\zeta\left(\frac{z'}{y_d} \right) \right]\frac{dz'}{|y_d|^{d-1}}\\
&=-\frac{1}{y_d} \int_{\mathbb{R}^{d-1}} (D_i \phi)(y'-y_d z')\left[(d-1) \zeta(z')+z'\cdot \nabla \zeta(z')\right] dz'.
\end{align*}
When $y_d =0$, it is easy to see that both quantities in \eqref{eq:yd-identity} become zero. Hence it follows that 
\[ |y_d D_{id} \phi^{(y_d)}(y')| \leq \theta \int_{\mathbb{R}^{d-1}} |(d-1)\zeta(z')+z' \cdot \nabla \zeta(z')| dz' = N\theta.\]
This proves the case $j=d$. If $j=1,\dots,d-1$, then 
\[
  |y_dD_{ij}\phi(y')|=\left|\int_{\mathbb{R}^{d-1}} (D_i \phi)(y'-y_d z') D_j \zeta(z')dz'\right|\leq N \theta
\]
for some constant $N=N(d,M)>0$.
Similarly, by using \eqref{eq:Dd-phi} and a change of variables,
\[ |y_d D_{dd} \phi^{(y_d)}(y')|\leq N \theta \]
for some constant $N=N(d,M)>0$. This proves (b).

(c) Observe that 
\begin{align*}
D_{ij}\phi^{(y_d)}(y')&=\int_{\mathbb{R}^{d-1}} (D_{ij} \phi)(y'-y_d z') \zeta(z')dz',\qquad i,j=1,\dots,d-1, \\
D_{id}\phi^{(y_d)}(y')&=-\int_{\mathbb{R}^{d-1}} (D_{il} \phi)(y'-y_d z') z_l' \zeta(z') dz',\qquad i=1,\dots,d-1, \\
D_{dd}\phi^{(y_d)}(y')&=\int_{\mathbb{R}^{d-1}} (D_{lk} \phi)(y'-y_d z') z_l'z_k' \zeta(z') dz',\qquad i=1,\dots, d-1.\label{eq:derivative-dd}
\end{align*}
Hence it follows that 
\[  |D^2 \phi^{(y_d)}(y')|\leq N.\]
Similarly, one can see that 
\[  |D^3 \phi^{(y_d)}(y')|\leq N\]
for some constant $N>0$, which proves (c).

(d) Since $\eta \in C_0^\infty(\mathbb{R}^{d-1})$ and $\phi \in C^{2,1}(\mathbb{R}^{d-1})$, the result follows from the identity \eqref{eq:yd-identity}. This completes the proof of Lemma \ref{lem:estimates-boundary}.
\end{proof}

By \eqref{eq:transformation-1}, we have 
\begin{equation*}\label{eq:Jacobian}
    (D\Phi)_{ij}(y)=\frac{\partial x_i}{\partial y_j}(y)=\left\{ 
    \begin{aligned}
        & \delta_{ij}-y_d D_{ij}\phi^{(y_d)}(y')&&\quad\text{if } i,j=1,\dots,d-1,\\
    & D_{j}\phi^{(y_d)}(y')&&\quad\text{if } i=d,\,j=1,\dots, d-1,\\
    & - D_i \phi^{(y_d)}(y')-y_d D_{id} \phi^{(y_d)}(y')&&\quad\text{if } i=1,\dots,d-1,\,j=d,\\
    & 1+D_d \phi^{(y_d)}(y')&&\quad\text{if } i=j=d,
    \end{aligned}\right.
\end{equation*}
and 
\[  \left(\frac{\partial x}{\partial y}\right)(0)=I.\]
 Hence it follows from the inverse function theorem that there exist $R^*>0$ and a $C^{2,1}$-diffeomorphism $\Psi:B_{R^*}(x_0)\rightarrow \Psi(B_{R^*}(x_0))$ whose inverse is $\Phi$. 
 
Define
\[  \hat{u}(t,y)=u(t,\Phi(y))\quad\text{and}\quad \tilde{u}(t,y)=\left(\frac{\partial x}{\partial y}\right)^T \hat{u}(t,y).\]

It is easy to check that when $y_d=0$,
\begin{equation}\label{eq:no-penetration}
   \tilde{u}(t,y)\cdot (0,\dots,0,-1)=0
\end{equation}
if and only if 
\[  \hat{u}(t,y)\cdot (\nabla \phi(y'),-1)=0 .  \]

Using this boundary flattening map, we connect the boundary condition \eqref{eq:general-Lions-boundary} with the Lions boundary conditions \eqref{eq:Lions-bdry-A}.

\begin{proposition}\label{prop:Lions-flatterning}
The boundary condition \eqref{eq:general-Lions-boundary} implies 
\begin{equation}\label{eq:Lions-bdry-A}
\tilde{u}^d(t,y',0)=0\quad \text{and}\quad    D_d \tilde{u}^k(t,y',0)=0
\end{equation}
for $(t,y') \in (-(R^*)^2,0]\times \Psi(\partial\Omega \cap B_{R^*}(x_0))$, where $k=1,\dots,d-1$. 
\end{proposition}
\begin{proof}
The chain rule gives
\begin{align*}
D_l\tilde{u}^k=\frac{\partial}{\partial y_l} \left(\frac{\partial x_i}{\partial y_k}\hat{u}^i\right)=\frac{\partial^2 x_i}{\partial y_k \partial y_l} \hat{u}^i + \frac{\partial x_i}{\partial y_k}\frac{\partial \hat{u}^i}{\partial y_l}
\end{align*}
and
\begin{align*}
D_k\tilde{u}^l=\frac{\partial}{\partial y_k} \left(\frac{\partial x_i}{\partial y_l}\hat{u}^i\right)=\frac{\partial^2 x_i}{\partial y_l \partial y_k} \hat{u}^i + \frac{\partial x_i}{\partial y_l}\frac{\partial \hat{u}^i}{\partial y_k},
\end{align*}
which imply that 
\begin{equation}\label{eq:Lions-calc-A}
D_l \tilde{u}^k(t,y)-D_k \tilde{u}^l(t,y)=\left(\frac{\partial x_i}{\partial y_k} \frac{\partial x_m}{\partial y_l}-\frac{\partial x_i}{\partial y_l} \frac{\partial x_m}{\partial y_k} \right)D_mu^i(t,\Phi(y)).
\end{equation}

Since $\Psi(\partial \Omega \cap B_{R^*}(x_0))\subset \{y_d =0\}$, we have
\begin{equation}\label{eq:Lions-calc-B}
\begin{aligned}
 \frac{\partial x_i}{\partial y_k}\frac{\partial x_m}{\partial y_d} D_m u^i= (-D_i \phi)D_i u^k+ (D_k \phi)(-D_i \phi) D_i u^d + D_d u^k + (D_k \phi)(D_d u^d)
\end{aligned}
\end{equation}
and 
\begin{equation}\label{eq:Lions-calc-C}
\begin{aligned}
 \frac{\partial x_i}{\partial y_d}\frac{\partial x_m}{\partial y_k} D_m u^i= (-D_i \phi)D_k u^i + D_k u^d - (D_i \phi)(D_k \phi)D_d u^i + (D_k\phi)(D_d u^d).
\end{aligned}
\end{equation}
Hence by \eqref{eq:Lions-calc-A}, \eqref{eq:Lions-calc-B}, and \eqref{eq:Lions-calc-C}, for $(t,y') \in (-(R^*)^2,0]\times \Psi(\partial\Omega \cap B_{R^*}(x_0))$ and $l=d$, we have
\begin{equation}\label{eq:flat-navier}
\begin{aligned}
D_d \tilde{u}^k-D_k\tilde{u}^d&=\sum_{i=1}^{d-1}(-D_i \phi)(D_i u^k-D_k u^i)+ D_d u^k - D_k u^d\\
&\relphantom{=} + \sum_{i=1}^{d-1}(D_k \phi)(-D_i \phi)(D_i u^d- D_d u^i)\\
&=\sum_{i=1}^{d-1} (-D_i \phi)[\omega_{ik}+(D_k \phi) \omega_{id}] + \omega_{dk},
\end{aligned}
\end{equation}
where $k=1,\dots,d-1$. Since 
\[ n_d=-\frac{1}{\sqrt{1+|\nabla \phi'(y')|^2}}\quad \text{ and }\quad n_j=\frac{D_j \phi(y')}{\sqrt{1+|\nabla \phi'(y')|^2}},\quad j=1,\dots, d-1,\]
it follows from the 2nd condition in \eqref{eq:general-Lions-boundary} that
\begin{equation}\label{eq:omega-condition}
\omega_{id}=\sum_{j=1}^{d-1} \omega_{ij} D_j \phi
\end{equation}
for $i=1,\dots,d$. Hence by \eqref{eq:flat-navier} and \eqref{eq:omega-condition}, we have
\begin{align*}
D_d\tilde{u}^k(t,y',0)-D_k\tilde{u}^d(t,y',0)&=\sum_{i=1}^{d-1} (-D_i \phi)\omega_{ik} + \omega_{dk}-(D_k\phi)\sum_{i=1}^{d-1} (D_i \phi)\omega_{id}\\
&=\omega_{kd}+\omega_{dk}-(D_k \phi) \sum_{i=1}^{d-1} (D_i \phi) \sum_{j=1}^{d-1} \omega_{ij} D_j \phi.
\end{align*}
Since $\omega_{ij}$ is antisymmetric and $(D_i \phi)(D_j\phi)$ is symmetric, it follows that 
\[ \sum_{i=1}^{d-1} (D_i \phi) \sum_{j=1}^{d-1} \omega_{ij} D_j \phi=0\quad \text{for $k=1,\dots,d-1$},\]
which implies
\[ D_d \tilde{u}^k(t,y',0)-D_k\tilde{u}^d(t,y',0)=0\]
for $(t,y')\in (-(R^*)^2,0]\times \Psi(\partial\Omega \cap B_{R^*}(x_0))$. Hence Proposition \ref{prop:Lions-flatterning} is proved.
\end{proof}

\subsection{Solvability results for Stokes equations with variable coefficients}\label{subsec:Stokes-variable}
We consider the following Stokes equations with variable coefficients
\begin{equation}\label{eq:Stokes-variable-half}
\left\{\begin{aligned}
\partial_t u -a^{ij}(t,x)D_{ij} u +\nabla p &=f&&\quad \text{in } (0,T)\times\mathbb{R}^d_+,\\
\Div u & =g&&\quad \text{in } (0,T)\times\mathbb{R}^d_+,
\end{aligned}
\right.
\end{equation}
subject to the Lions boundary conditions
\begin{equation}\label{eq:Lions}
u^d=0\quad \text{and}\quad D_d u^k=0,\quad 1\leq k\leq d-1
\end{equation}
 on $(0,T)\times\mathbb{R}^{d-1}\times\{0\}$. For simplicity, we write $\Omega_T=(0,T)\times\mathbb{R}^d_+$.

To state the solvability result, for a function $g$ and a vector field $G$, we say that $g_t = \Div G$  in $U\subset \Omega_T$ in the sense of distribution if 
\begin{equation}\label{eq:compatibility-condition}
 \int_{U} g \partial_t \phi \myd{x}dt = \int_{U} G \cdot D\phi \myd{x}dt 
 \end{equation}
holds for all $\phi \in C_0^\infty(U)$.

The following result can be found in \cite[Theorem 8.1]{DK23}. We refer the reader to \cite[Section 2]{DK23} for  the definitions of $\oSob{1,2}{s,q,w}(\Omega_T)$ and $\mathring{\mathcal{H}}^1_{s,q,w}(\Omega_T)$. 

\begin{theorem}\label{thm:VMO-half}
Let $T\in (0,\infty)$, $s,q\in (1,\infty)$, $K_0\geq 1$, $w=w_1(x)w_2(t)$, where
\[  w_1 \in A_q(\mathbb{R}^d,dx),\quad w_2 \in A_s(\mathbb{R},dt),\quad [w]_{A_{s,q}}\leq K_0.\]
There exists a constant $\gamma=\gamma(d,\nu,s,q,K_0)\in (0,1)$ such that under Assumption \ref{assump:VMO} $(\gamma)$, if $(u,p)$ is a strong solution to \eqref{eq:Stokes-variable-half} in $\Omega_T$ subject to the Lions boundary conditions \eqref{eq:Lions} and $u(0,\cdot)=0$ on $\mathbb{R}^d_+$ satisfying $(u,\nabla p)\in \oSob{1,2}{s,q,w}(\Omega_T)^d\times \Leb{s,q,w}(\Omega_T)^d$, $f\in \Leb{s,q,w}(\Omega_T)^d$,  $g\in \mathring{\mathcal{H}}^{1}_{s,q,w}(\Omega_T)$ satisfying $g_t=\Div G$  in $\Omega_T$ for some $G\in\Leb{s,q,w}(\Omega_T)^d$, then there exist constants $N_1=N_1(d,\nu,s,q,K_0)>0$ and $N_2=N_2(d,\nu,s,q,K_0,R_0)>0$ such that 
\begin{equation}\label{eq:a-priori-estimate-half}
\begin{aligned}
&\norm{\partial_t u}{\Leb{s,q,w}(\Omega_T)}+\norm{D^2 u}{\Leb{s,q,w}(\Omega_T)}+\norm{\nabla p}{\Leb{s,q,w}(\Omega_T)}\\
&\leq N_1\left(\norm{f}{\Leb{s,q,w}(\Omega_T)}+\norm{Dg}{\Leb{s,q,w}(\Omega_T)}+\norm{G}{\Leb{s,q,w}(\Omega_T)} \right)+N_2\norm{u}{\Leb{s,q,w}(\Omega_T)}.
\end{aligned}
\end{equation}
Furthermore, for every $f\in \Leb{s,q,w}(\Omega_T)^d$,  $g\in \mathring{\mathcal{H}}^{1}_{s,q,w}(\Omega_T)$ satisfying $g_t=\Div G$ in $\Omega_T$ for some $G\in\Leb{s,q,w}(\Omega_T)^d$ in the sense of \eqref{eq:compatibility-condition}, then there exists a unique solution $(u, p) \in \oSob{1,2}{s,q,w}(\Omega_T)^d\times \Leb{s,q,w,\mathrm{loc}}(\Omega_T)\slash\mathbb{R}$ of \eqref{eq:Stokes-variable-half} in $\Omega_T$ satisfying $\nabla p\in \Leb{s,q,w}(\Omega_T)$, the Lions boundary conditions \eqref{eq:Lions}, and $u(0,\cdot)=0$ on $\mathbb{R}^d_+$. Moreover, there exists a constant $N_3=N_3(d,\nu,s,q,K_0,R_0,T)>0$ such that
\begin{equation*}
\begin{aligned}
& \norm{u}{\Sob{1,2}{s,q,w}(\Omega_T)}+\norm{\nabla p}{\Leb{s,q,w}(\Omega_T)}\\
&\leq N_3 \left(\norm{f}{\Leb{s,q,w}(\Omega_T)}+\norm{Dg}{\Leb{s,q,w}(\Omega_T)}+\norm{G}{\Leb{s,q,w}(\Omega_T)}\right).
\end{aligned}
\end{equation*}
\end{theorem}

\subsection{Hessian estimates for Stokes equations near the flat boundary}
In this subsection, we recall some known results on Hessian estimates for Stokes equations in nondivergence form:
\begin{equation}\label{eq:flat-space}
\partial_t u -a^{ij}D_{ij} u + \nabla p =f\quad\text{and}\quad \Div u=g\quad \text{in } U,
\end{equation}
where $U$ is a cylindrical domain in $\mathbb{R}^{d+1}$.

The following interior regularity estimate was shown in Dong-Kwon \cite[Theorem 2.7]{DK23}. See Section \ref{sec:perturbed-Stokes} for the proof which contains the idea of the proof by using Dong-Phan \cite[Lemma 4.13]{DP21}.

\begin{theorem}\label{thm:interior-estimate-Stokes}
Let $s,q,q_0\in (1,\infty)$ and $r\in (0,R)$. There exists a constant $\gamma\in (0,1)$, depending on $d$, $\nu$, $s$, $q$, and $q_0)$, such that under Assumptions \ref{assump:VMO} $(\gamma)$, if $(u,p)\in \Sob{1,2}{q_0}(Q_{R})^d\times \Sob{0,1}{1}(Q_{R})$ is a solution of \eqref{eq:flat-space} in $Q_{R}$ for some $f \in \Leb{s,q}(Q_{R})^d$ and $g\in \Sob{0,1}{s,q}(Q_{R})$, then $D^2 u \in \Leb{s,q}(Q_r)$. Moreover, there exists a constant $N=N(d,s,q,\nu,r,R,R_0)>0$ such that 
\[ \norm{D^2 u}{\Leb{s,q}(Q_r)}\leq N\left( \norm{u}{\Leb{s,1}(Q_R)}+\norm{f}{\Leb{s,q}(Q_R)}+\norm{Dg}{\Leb{s,q}(Q_R)}\right).\]
\end{theorem}

For the boundary estimate, we use the following theorem which was proved in Dong-Kim-Phan \cite[Lemma 5.5]{DKP22} when $R\in [1/2,1]$. The case when $R\in [R_1/2,1]$ can be proved by slightly modifying the proof there.
\begin{theorem}\label{thm:iteration-preparation}
Let $R_1 \in (0,R_0)$, $R\in [R_1/2,1]$, $\gamma \in (0,1)$, $\kappa \in (0,1/4)$, $s,q\in (1,\infty)$, $q_1 \in (1,\min\{s,q\})$, and $q_0 \in (1,q_1)$. Let $(u,p)\in \Sob{1,2}{s,q}(Q_{R+R_1}^+)^d \times \Sob{0,1}{1}(Q_{R+R_1}^+)$ be a strong solution to \eqref{eq:flat-space} with the boundary conditions \eqref{eq:Navier-flat} on $Q_{R+R_1} \cap \{(t,x) : x_d=0\}$, where $f\in \Leb{s,q}(Q_{R+R_1}^+)^d$, and $Dg\in \Leb{s,q}(Q_{R+R_1}^+)^d$. Then under Assumption \ref{assump:VMO} $(\gamma)$, we have 
\begin{align*}
\norm{D^2 u}{\Leb{s,q}(Q_R^+)}&\leq N\kappa^{-(d+2)/q_0}\left[ \norm{f}{\Leb{s,q}(Q_{R+R_1/2}^+)}+\norm{Dg}{\Leb{s,q}(Q_{R+R_1/2}^+)}\right.\\
&\relphantom{=}\left.+R_1^{-1} \norm{Du}{\Leb{s,q}(Q_{R+R_1/2}^+)}\right]\\
&\relphantom{=}+N\left(\kappa^{-(d+2)/q_0}\gamma^{1/q_0-1/q_1}+\kappa^{1/2}\right)\norm{D^2u}{\Leb{s,q}(Q_{R+R_1/2}^+)}
\end{align*}
for some constant $N=N(d,s,q,\nu)>0$.
\end{theorem}

\section{Perturbed Stokes equations}\label{sec:perturbed-Stokes}
We consider the following perturbed Stokes systems:
\begin{equation}\label{eq:perturbed-Stokes-I}
\left\{
\begin{alignedat}{2}
\partial_t u -a^{ij}(t,x)D_{ij} u +b^i(t,x)D_i u + c(t,x)u+\nabla p &=f&&\quad \text{in }U,\\
\Div (\widehat{A}u)&=g&&\quad \text{in } U,\\ 
\end{alignedat}
\right.
\end{equation}
where $U=(S,T)\times \Omega$ is a cylindrical domain in $\mathbb{R}^{d+1}$. Here we assume that $a^{ij}:\mathbb{R}\times \Omega\rightarrow\mathbb{R}^{d \times d}$ is uniformly elliptic \eqref{eq:elliptic}, $b:\mathbb{R}\times \Omega\rightarrow\mathbb{R}^d$, and $c:\mathbb{R}\times \Omega\rightarrow\mathbb{R}$ are bounded {by a constant $K$}. The matrix $\widehat{A}:\Omega \rightarrow \mathbb{R}^{d\times d}$ is symmetric and satisfies the following regularity assumption.

\begin{assumption}[$\beta$]\label{assump:delta-N}
There exist constants $R_0,N>0$ such that 
\[
 \norm{\widehat{A}-I}{\Leb{\infty}(B_{R_0}^+)}\leq \beta\quad \text{ and }\quad   \sum_{k=1}^2\norm{D^k \widehat{A}}{\Leb{\infty}(B_{R_0}^+)}\leq N.
\]
\end{assumption}

The following theorem is the main result of this section.
\begin{theorem}\label{thm:main-theorem-perturb}
Let $s,q,q_0 \in (1,\infty)$ and $R>0$. There exist constants $\beta,\gamma>0$ depending on $d,\nu,s,q$ such that under Assumptions \ref{assump:VMO} $(\gamma)$ and \ref{assump:delta-N} $(\beta)$, if $(u,p)\in \Sob{1,2}{q_0}(Q_{R}^+)^d\times \Sob{0,1}{1}(Q_{R}^+)$ is a solution of \eqref{eq:perturbed-Stokes-I} in $Q_{R}^+$  with the Lions boundary conditions \eqref{eq:Lions} on $(-R^2,0]\times B_{R}' \times\{0\}$ for some $f \in \Leb{s,q}(Q_{R}^+)^d$ and $g\in \Sob{0,1}{s,q}(Q_{R}^+)$, then $D^2 u \in \Leb{s,q}(Q_{R/2}^+)$. Moreover, there exists a constant $N=N(d,s,q,\nu,R,R_0)>0$ such that 
\begin{equation*}
\norm{D^2 u}{\Leb{s,q}(Q_{R/2}^+)}\leq N\left(\norm{u}{\Leb{s,1}(Q_R^+)}+\norm{f}{\Leb{s,q}(Q_R^+)}+\norm{Dg}{\Leb{s,q}(Q_R^+)}\right).
\end{equation*}
\end{theorem}
 
\begin{remark}\label{rem:u-t-estimate}
In the previous paper \cite{DK23}, we assumed that $u,Du,D^2u\in\Leb{q_0}(Q_R^+)$ for some $q_0>1$ and $\partial_t u \in \Leb{1}(Q_R^+)$. However, since we will use the maximal regularity estimates for \eqref{eq:perturbed-Stokes-I}, we need to assume that $\partial_t u \in \Leb{q_0}(Q_R^+)$. See \eqref{eq:perturbed-Stokes-II} and \eqref{eq:u-1-epsilon-perturb}. Also, from the equation, one can see that $Dp \in \Leb{q_0}(Q_{R}^+)^d$.
\end{remark}

In Section \ref{sub:perturb-1}, we first study the weighted solvability of the Cauchy problem of \eqref{eq:perturbed-Stokes-I} on $\Omega_T=(0,T)\times \mathbb{R}^d_+$ which is necessary to obtain regularity improving lemma (Lemma \ref{lem:regularity-lemma-Stokes}) that will be used in the proof of Theorem \ref{thm:main-theorem-perturb}. Next, we obtain a \emph{priori local boundary Hessian estimates} for perturbed Stokes equations \eqref{eq:perturbed-Stokes-I} using Theorem \ref{thm:iteration-preparation}. The proof of Theorem \ref{thm:main-theorem-perturb} will be given in Section \ref{subsec:perturb-3}.

\subsection{Weighted solvability results for perturbed Stokes systems}\label{sub:perturb-1} In this subsection, we prove a unique solvability result for \eqref{eq:perturbed-Stokes-I} on $\Omega_T=(0,T)\times\mathbb{R}^d_+$ in weighted Sobolev spaces $\Sob{1,2}{s,q,w}(\Omega_T)$ with the Lions boundary conditions \eqref{eq:Lions}. We impose the following assumption on $\widehat{A}:\mathbb{R}^d \rightarrow\mathbb{R}^{d\times d}$.
\begin{assumption}[$\beta$]\label{assump:delta-N-1}
There exists a constant $N>0$ such that  
\[  \norm{\widehat{A}-I}{\Leb{\infty}(\mathbb{R}^d)}\leq \beta\quad\text{and}\quad \sum_{k=1}^2\norm{D^k \widehat{A}}{\Leb{\infty}(\mathbb{R}^d)}\leq N.\]
\end{assumption}

\begin{theorem}\label{thm:solvability-modified-divergence}
Let $T\in (0,\infty)$, $s,q\in (1,\infty)$, $K_0\geq 1$, $w=w_1(x)w_2(t)$, where
\[  w_1 \in A_q(\mathbb{R}^d,dx),\quad w_2 \in A_s(\mathbb{R},dt),\quad [w]_{A_{s,q}}\leq K_0.\]
There exist $\gamma=\gamma(d,\nu,s,q,K_0)$ and $\beta = \beta(d,\nu,s,q,K_0)$ such that under Assumptions \ref{assump:VMO} $(\gamma)$ and \ref{assump:delta-N-1} $(\beta)$, for any $f \in \Leb{s,q,w}(\Omega_T)^d$ and $g\in \mathring{\mathcal{H}}^{1}_{s,q,w}(\Omega_T)$ satisfying the compatibility condition $g_t = \Div G$ in $\Omega_T$ in the sense of distribution for some $G\in \Leb{s,q,w}(\Omega_T)^d$, there exists a unique solution $(u,p)\in\oSob{1,2}{s,q,w}(\Omega_T)^d\times \Leb{s,q,w,\mathrm{loc}}(\Omega_T)\slash\mathbb{R}$ satisfying $\nabla p\in \Leb{s,q,w}(\Omega_T)^d$, \eqref{eq:perturbed-Stokes-I}, and the Lions boundary conditions \eqref{eq:Lions}. Moreover, we have 
\begin{equation}\label{eq:a-priori-estimate-half-2}
\begin{aligned}
&\norm{u}{\Sob{1,2}{s,q,w}(\Omega_T)}+\norm{\nabla p}{\Leb{s,q,w}(\Omega_T)}\\
&\leq N\left(\norm{f}{\Leb{s,q,w}(\Omega_T)}+\norm{Dg}{\Leb{s,q,w}(\Omega_T)}+\norm{G}{\Leb{s,q,w}(\Omega_T)}\right)
\end{aligned}
\end{equation}
for some constant $N=N(d,s,q,\nu,K_0,R_0,T)>0$.
\end{theorem}
\begin{proof}
By the method of continuity together with Theorem \ref{thm:VMO-half}, it suffices to obtain a priori estimate \eqref{eq:a-priori-estimate-half-2}. We first note that if $u\in \Sob{1,2}{s,q,w}(\Omega_T)^d$ satisfies \eqref{eq:Lions}, then for each $\varepsilon \in (0,1)$, we have
\begin{equation}\label{eq:half-extension-weight}
\norm{Du}{\Leb{s,q,w}(\Omega_T)}\leq \varepsilon\norm{D^2u}{\Leb{s,q,w}(\Omega_T)}+\frac{N}{\varepsilon}\norm{u}{\Leb{s,q,w}(\Omega_T)}
\end{equation}
for some constant $N=N(d,s,q,K_0)>0$.  Indeed, let $\tilde{u}^k$ be the even extensions of $u^k$ with respect to $x_d$ for $k=1,\dots,d-1$, and let $\tilde{u}^d$ be the odd extension of $u^d$ with respect to $x_d$. Then one can easy to see that $\tilde{u}\in \Sob{1,2}{s,q,w}((0,T)\times\mathbb{R}^d)^d$. Hence we can apply \eqref{eq:weighted-interpolation} to $\tilde{u}$ to get \eqref{eq:half-extension-weight}.

If we write $h=\Div ((I-\widehat{A})u)+g$, then it is easy to see that 
\[  \int_{\Omega_T} h \varphi_t dxdt = \int_{\Omega_T} \left[(I-\widehat{A}) u_t +G\right] \cdot \nabla \varphi \myd{x}dt \]
for all $\varphi \in C_0^\infty((0,T)\times\mathbb{R}^d_+)$. Since $\Div u =h$ and 
\[ \partial_t u -a^{ij}D_{ij}u+\nabla p =f+b^i D_i u +cu,\]
 it follows from \eqref{eq:a-priori-estimate-half} that 
\begin{align*}
&\norm{\partial_t u}{\Leb{s,q,w}(\Omega_T)}+\norm{D^2 u}{\Leb{s,q,w}(\Omega_T)}+\norm{\nabla p}{\Leb{s,q,w}(\Omega_T)}\\
&\leq N_1\left( \norm{f}{\Leb{s,q,w}(\Omega_T)}+\norm{D^2((I-\widehat{A})u)}{\Leb{s,q,w}(\Omega_T)}+\norm{(\widehat{A}-I)\partial_tu}{\Leb{s,q,w}(\Omega_T)}\right.\\
&\relphantom{=}\left.+\norm{Dg}{\Leb{s,q,w}(\Omega_T)}+\norm{G}{\Leb{s,q,w}(\Omega_T)}+\norm{Du}{\Leb{s,q,w}(\Omega_T)}\right)+N_2\norm{u}{\Leb{s,q,w}(\Omega_T)}\\
&\leq N_1\left( \norm{f}{\Leb{s,q,w}(\Omega_T)}+\norm{Du}{\Leb{s,q,w}(\Omega_T)}+\beta\left(\norm{\partial_t u}{\Leb{s,q,w}(\Omega_T)}+\norm{D^2u}{\Leb{s,q,w}(\Omega_T)}\right)\right.\\
&\relphantom{=}+\left.\norm{Dg}{\Leb{s,q,w}(\Omega_T)}+\norm{G}{\Leb{s,q,w}(\Omega_T)}\right)+N_2\norm{u}{\Leb{s,q,w}(\Omega_T)},
\end{align*}
where the constants $N_1=N_1(d,\nu,s,q,K_0)>0$ and $N_2=N_2(d,\nu,s,q,K_0,R_0)>0$. Hence by the interpolation inequality \eqref{eq:half-extension-weight}, we have
\begin{align*}
&\norm{\partial_t u}{\Leb{s,q,w}(\Omega_T)}+\norm{D^2 u}{\Leb{s,q,w}(\Omega_T)}+\norm{\nabla p}{\Leb{s,q,w}(\Omega_T)}\\
&\leq N_1 \norm{f}{\Leb{s,q,w}(\Omega_T)}+\frac{1}{4}\norm{D^2u}{\Leb{s,q,w}(\Omega_T)}+N_1\beta\left(\norm{\partial_t u}{\Leb{s,q,w}(\Omega_T)}+\norm{D^2u}{\Leb{s,q,w}(\Omega_T)}\right)\\
&\relphantom{=}+N_2\norm{u}{\Leb{s,q,w}(\Omega_T)}+N_1\norm{Dg}{\Leb{s,q,w}(\Omega_T)}+N_1\norm{G}{\Leb{s,q,w}(\Omega_T)}.
\end{align*}
 Now choose $\beta \in (0,1)$ so that $N_1\beta \leq 1/4$. Then we get 
\begin{align*}
&\norm{\partial_t u}{\Leb{s,q,w}(\Omega_T)}+\norm{D^2 u}{\Leb{s,q,w}(\Omega_T)}+\norm{\nabla p}{\Leb{s,q,w}(\Omega_T)}\\
&\leq N\left(\norm{f}{\Leb{s,q,w}(\Omega_T)}+\norm{u}{\Leb{s,q,w}(\Omega_T)}+\norm{Dg}{\Leb{s,q,w}(\Omega_T)}+\norm{G}{\Leb{s,q,w}(\Omega_T)}\right),
\end{align*}
where $N=N(d,s,q,\nu,K_0,R_0)>0$. Hence by using the time splitting argument as in the proof of \cite[Theorem 2.5]{DK23}, we get 
\[
\norm{u}{\Sob{1,2}{s,q,w}(\Omega_T)}+\norm{\nabla p}{\Leb{s,q,w}(\Omega_T)}\leq N\left(\norm{f}{\Leb{s,q,w}(\Omega_T)}+\norm{Dg}{\Leb{s,q,w}(\Omega_T)}+\norm{G}{\Leb{s,q,w}(\Omega_T)}\right)
\]
for some constant $N=N(d,s,q,\nu,K_0,R_0,T)>0$. This completes the proof of Theorem \ref{thm:solvability-modified-divergence}.
\end{proof}

As an application, we have the following regularity improving lemma. A simplified version of this lemma has already been proven in \cite[Lemma 7.3]{DK23}. We omit its proof since the argument is similar. 
\begin{lemma}\label{lem:regularity-lemma-Stokes}
Let $T\in (0,\infty)$ and $s,q,q_0\in (1,\infty)$. There exist constants $\beta,\gamma>0$ depending on $d,\nu,s,q_0,q$ such that under Assumptions \ref{assump:VMO} $(\gamma)$ and \ref{assump:delta-N-1} $(\beta)$, if $(u,p)$ is a strong solution to \eqref{eq:perturbed-Stokes-I} in $\Omega_T$ with $u(0,\cdot)=0$ on $\mathbb{R}^d_+$ with the Lions boundary conditions \eqref{eq:Lions} satisfying
\[
u \in \oSob{1,2}{q_0}(\Omega_T)^d,\quad \nabla p\in \Leb{q_0}(\Omega_T)^d,\quad
\]$f\in \Leb{\infty}(\Omega_T)^d$ having compact support in $\Omega_T${,} and a bounded $g$  satisfying $g_t=\Div G$ in $\Omega_T$, where $G \in \Sob{0,2}{\infty}(\Omega_T)$, and both $g$ and $G$ has compact support in $\Omega_T$, then $(u,\nabla p)\in \Sob{1,2}{s,q}(\Omega_T)^d \times \Leb{s,q}(\Omega_T)^d$.
\end{lemma}

\subsection{A priori regularity estimates for perturbed Stokes systems}\label{subsec:sub:perturb-2}
Next, we obtain a priori local boundary mixed norm Hessian estimates for \eqref{eq:perturbed-Stokes-I} under the Lions boundary conditions. 

\begin{theorem}\label{thm:regularity-a-priori-local}
Let $s,q\in(1,\infty)$ and $R>0$. There exist $\beta,\gamma>0$ depending on $d,\nu,s,q$ such that under Assumptions \ref{assump:VMO} $(\gamma)$ and \ref{assump:delta-N} $(\beta)$, if $(u,p)\in \tilde{W}^{1,2}_{s,q}(Q_{R}^+)^d\times \Sob{0,1}{1}(Q_{R}^+)$ is a solution of \eqref{eq:perturbed-Stokes-I} in $Q_{R}^+$  under the Lions boundary conditions \eqref{eq:Lions} on $(-R^2,0]\times B_{R}' \times\{0\}$ for some $f \in \Leb{s,q}(Q_{R}^+)^d$ and $g\in \Sob{0,1}{s,q}(Q_{R}^+)$, then 
\begin{align*}
 \norm{D^2 u}{\Leb{s,q}(Q_{3R/4}^+)}&\leq N\left(\norm{u}{\Leb{s,1}(Q_R^+)}+\norm{f}{\Leb{s,q}(Q_R^+)}+\norm{Dg}{\Leb{s,q}(Q_R^+)}\right)
\end{align*}
holds for some constant $N=N(d,s,q,\nu,R,R_0)>0$.
\end{theorem}
\begin{proof}
By a partition of unity argument, we only prove the case when $0<R\leq R_0$, where $R_0$ is the number in Assumptions \ref{assump:VMO} $(\gamma)$ and  \ref{assump:delta-N} $(\beta)$. Define
\[  r_k=R-\frac{R}{2^{k+1}},\quad k=1,2,\dots.\]
Then $r_1=3R/4$ and $r_k$ is increasing satisfying $\lim_{k\rightarrow\infty} r_k=R$. Since $(u,p)$ satisfies
\begin{equation}\label{eq:Stokes-nondiv-rescaled}
\left\{\begin{aligned}
\partial_t u-a^{ij}D_{ij} u +\nabla p&=\tilde{f}\\
\Div u &=\tilde{g},
\end{aligned}
\right.
\end{equation}
where
\[
\tilde{f}=f+b^iD_i u +cu\quad\text{and}\quad \tilde{g}=g+\Div((I-\widehat{A})u),
\]
we apply Theorem \ref{thm:iteration-preparation} to \eqref{eq:Stokes-nondiv-rescaled} with $R=r_k$ and $R_1=2^{-k-2}R$. Since $r_k+R_1=r_{k+1}$, we have
\begin{equation}\label{eq:Hessian-u}
\begin{aligned}
\norm{D^2 u}{\Leb{s,q}(Q_{r_k}^+)}&\leq N \kappa^{-(d+2)/q_0} (\norm{\tilde{f}}{\Leb{s,q}(Q_{r_{k+1}}^+)} +  \norm{D\tilde{g}}{\Leb{s,q}(Q_{r_{k+1}}^+)})\\
&\relphantom{=}+N\left(\kappa^{-(d+2)/q_0} \gamma^{1/q_0-1/q_1}+\kappa^{1/2} \right)\norm{D^2 u}{\Leb{s,q}(Q_{r_{k+1}}^+)}\\
&\relphantom{=}+N\kappa^{-(d+2)/q_0} \frac{2^{k}}{R} \norm{Du}{\Leb{s,q}(Q_{r_{k+1}}^+)}
\end{aligned}
\end{equation}
for some constant $N=N(d,s,q,\nu)>0$.

Note that 
\begin{equation}\label{eq:D-tilde-f}
\begin{aligned}
\norm{\tilde{f}}{\Leb{s,q}(Q_{r_{k+1}}^+)}&\leq N \left(\norm{f}{\Leb{s,q}(Q_{r_{k+1}}^+)}+\norm{Du}{\Leb{s,q}(Q_{r_{k+1}}^+)}+\norm{u}{\Leb{s,q}(Q_{r_{k+1}}^+)}\right)
\end{aligned}
\end{equation}
for some constant $N=N(d,s,q,K,\nu)>0$. 
 
Since $\widehat{A}$ satisfies Assumption \ref{assump:delta-N} $(\beta)$, it follows that \begin{equation}\label{eq:D-tilde-g}
\begin{aligned}
\norm{D\tilde{g}}{\Leb{s,q}(Q_{r_{k+1}}^+)}&\leq \norm{Dg}{\Leb{s,q}(Q_{r_{k+1}}^+)}+\norm{D\Div((I-\widehat{A})u)}{\Leb{s,q}(Q_{r_{k+1}}^+)}\\
&\leq \norm{Dg}{\Leb{s,q}(Q_{r_{k+1}}^+)}+N(\norm{Du}{\Leb{s,q}(Q_{r_{k+1}}^+)}+\norm{u}{\Leb{s,q}(Q_{r_{k+1}}^+)})\\
&\relphantom{=}+N\beta\norm{D^2 u}{\Leb{s,q}(Q_{r_{k+1}}^+)},
\end{aligned}
\end{equation}
where $N=N(d,s,q,\nu)>0$. 

On the other hand, by the Poincar\'e inequality, we have 
\[
  \norm{u}{\Leb{s,q}(Q_{r_{k+1}}^+)}\leq N\left(\norm{Du}{\Leb{s,q}(Q_{r_{k+1}}^+)} + r_{k+1}^{d(1/q-1)}\norm{u}{\Leb{s,1}(Q_{r_{k+1}}^+)} \right)
\]
for some constant $N=N(d,q)>0$. Hence it follows from \eqref{eq:Hessian-u}, \eqref{eq:D-tilde-f}, and \eqref{eq:D-tilde-g} that  
\begin{equation}\label{eq:Hessian-ak-before}
\begin{aligned}
\norm{D^2 u}{\Leb{s,q}(Q_{r_k}^+)}&\leq N \left(\kappa^{1/2}+\kappa^{-(d+2)/q_0}(\gamma^{1/q_0-1/q_1}+\beta) \right)\norm{D^2 u}{\Leb{s,q}(Q_{r_{k+1}}^+)}\\
&\relphantom{=}+N\kappa^{-(d+2)/q_0}\left(\norm{f}{\Leb{s,q}(Q_{r_{k+1}}^+)}+\norm{Dg}{\Leb{s,q}(Q_{r_{k+1}}^+)}\right)\\
&\relphantom{=}+N\kappa^{-(d+2)/q_0}\left(1+\frac{2^k}{R} \right)\norm{Du}{\Leb{s,q}(Q_{r_{k+1}}^+)}\\
&\relphantom{=}+N\kappa^{-(d+2)/q_0} r^{d(1/q-1)}_{k+1}\norm{u}{\Leb{s,1}(Q_{r_{k+1}}^+)}\end{aligned}
\end{equation}
for some constant $N=N(d,s,q,\nu)>0$. 

Moreover, by the interpolation inequality \eqref{eq:Du-interpolation-L1}, for each $\varepsilon\in(0,1)$, we have
\begin{equation}\label{eq:Hessian-ak}
\begin{aligned}
&\kappa^{-(d+2)/q_0} a_k\norm{Du}{\Leb{s,q}(Q_{r_{k+1}}^+)}\leq N\varepsilon \norm{D^2u}{\Leb{s,q}(Q_{r_{k+1}}^+)}\\
&\relphantom{=}+N\left(r^{-1+d(1/q-1)}_{k+1}\kappa^{-(d+2)/q_0}a_k  + \left(\kappa^{-\frac{d+2}{q_0}} a_k\right)^{\frac{1}{1-\theta}} \varepsilon^{-\frac{\theta}{1-\theta}} \right)\norm{u}{\Leb{s,1}(Q_{r_{k+1}}^+)},
\end{aligned}
\end{equation}
where $N=N(d,s,q,\nu)>0$,
\[ a_k=1+\frac{2^k}{R},\quad \text{and}\quad \theta=\frac{1+1/d-1/q}{1+2/d-1/q}\in [1/2,1).\]
 By combining \eqref{eq:Hessian-ak-before} with \eqref{eq:Hessian-ak}, we get 
\begin{equation}\label{eq:Hessian-ak-later}
\begin{aligned}
&\norm{D^2 u}{\Leb{s,q}(Q_{r_k}^+)}\leq N_1(\varepsilon+\kappa^{1/2}+\kappa^{-(d+2)/q_0}(\gamma^{1/q_0-1/q_1}+\beta) )\norm{D^2u}{\Leb{s,q}(Q_{r_{k+1}}^+)}\\
&\relphantom{=}+N_1\kappa^{-(d+2)/q_0}\left(\norm{f}{\Leb{s,q}(Q_{r_{k+1}}^+)}+\norm{Dg}{\Leb{s,q}(Q_{r_{k+1}}^+)}\right)\\
&\relphantom{=}+N_2\left(\kappa^{-(d+2)/q_0}a_k +\left(\kappa^{-\frac{d+2}{q_0}} a_k\right)^{\frac{1}{1-\theta}} \varepsilon^{-\frac{\theta}{1-\theta}} \right)\norm{u}{\Leb{s,1}(Q_{r_{k+1}}^+)},
\end{aligned}
\end{equation}
where $N_1=N_1(d,s,q,\nu)>0$ and $N_2=N_2(d,s,q,\nu,R)>0$.

Choose $\kappa$ sufficiently small so that $N\kappa^{1/2} \leq 1/36$ and then choose $\gamma$ and $\beta$ sufficiently small so that 
$$
N\kappa^{-(d+2)/q_0} \gamma^{1/q_0-1/q_1}\leq 1/36,\quad 
N\kappa^{-(d+2)/q_0}\beta \leq 1/36.
$$ 
Finally, choose $\varepsilon>0$ so that $N\varepsilon \leq 1/36$. If we set 
\[ \mathsf{A}_k=\left(\frac{1}{9}\right)^k \norm{D^2 u}{\Leb{s,q}(Q_{r_k}^+)},\]
then by multiplying $9^{-k}$ to the inequality \eqref{eq:Hessian-ak-later} and taking summation over $k=1$, $2$, ..., we get 
\begin{align*}
\sum_{k=1}^\infty \mathsf{A}_k &\leq \sum_{k=1}^\infty \mathsf{A}_{k+1}+N_1 \left(\norm{f}{\Leb{s,q}(Q_R^+)}+ \norm{Dg}{\Leb{s,q}(Q_R^+)}\right)\\
&\relphantom{=}+N_2\left(1+\frac{1}{R}+\frac{1}{R^{1/(1-\theta)}}\right)\norm{u}{\Leb{s,1}(Q_R^+)}
\end{align*}
for some constants $N_1=N_1(d,s,q,\nu)>0$ and $N_2=N_2(d,s,q,\nu,R)>0$. Hence by cancelling the summation on the right-hand side, we get 
\begin{align*}
 \norm{D^2 u}{\Leb{s,q}(Q_{3R/4}^+)}&\leq N\left(\norm{u}{\Leb{s,1}(Q_R^+)}+\norm{f}{\Leb{s,q}(Q_R^+)}+\norm{Dg}{\Leb{s,q}(Q_R^+)}\right),
\end{align*}
where $N=N(d,s,q,\nu,R)>0$. This completes the proof of Theorem \ref{thm:regularity-a-priori-local}.
\end{proof}

\subsection{Proof of Theorem \ref{thm:main-theorem-perturb}}\label{subsec:perturb-3} This subsection is devoted to the proof of Theorem \ref{thm:main-theorem-perturb}. 
\begin{proof}[Proof of Theorem \ref{thm:main-theorem-perturb}]
Suppose that $(u,p) \in \Sob{1,2}{q_0}(Q_{R}^+)^d \times \Sob{0,1}{1}(Q_{R}^+)$ satisfies the Lions boundary condition. Let $\tilde{u}^k$ be the even extension of $u^k$ with respect to $x_d$, $k=1,\dots,d-1$, $\tilde{u}^d$ be the odd extension of $u^d$ with respect to $x_d$. Let $\tilde{f}^k(t,\cdot)$ be the even extension of $f^k(t,\cdot)$ for $k=1,\dots,d-1$, and $\tilde{f}^d(t,\cdot)$ be the odd extension of $f^d(t,\cdot)$. Similarly, let $\tilde{g}(t,\cdot)$ be the even extension of $g(t,\cdot)$ with respect to $x_d$. Let $\tilde{p}$ be the even extension of $p$ in $x_d$. Then $\tilde{u}\in \Sob{1,2}{q_0}(Q_{R})^d$, $\tilde{p} \in \Sob{0,1}{1}(Q_{R})$, $\tilde{f}\in \Leb{s,q}(Q_{R})^d$, and $\tilde{g}\in \tilde{W}^{0,1}_{s,q}(Q_{R})$. 

Define $\overline{a}^{ij}(t,x',x_d)=a^{ij}(t,x',x_d)$ if $x_d>0$. For $x_d<0$, define 
\[ 
\overline{a}^{ij}(t,x',x_d)=\begin{cases}
a^{ij}(t,x',-x_d)& \text{ if $i,j=1,\dots,d-1$,}\\
-a^{id}(t,x',-x_d)& \text{ if $i=1,\dots,d-1$,  $j=d$,}\\
-a^{dj}(t,x',-x_d)& \text{ if $i=d$,  $j=1,\dots,d-1$,}\\
a^{dd}(t,x',-x_d) & \text{ if $i=j=d$.}
\end{cases}
\] 
Similarly, we define $\overline{b}^i(t,x',x_d)=b^i(t,x',x_d)$ and $\overline{c}(t,x',x_d)=c(t,x',x_d)$ if $x_d>0$ and 
\[ \overline{b}^i(t,x',x_d)=\left\{
\begin{aligned}
b^i(t,x',-x_d) & \text{ if $i=1,\dots,d-1$}\\
-b^d(t,x',-x_d) & \text{ if $i=d$}    
\end{aligned}\right.,
\quad 
\overline{c}(t,x',x_d)=c(t,x',-x_d).
\]
Also, let $\overline{\widehat{A}}^{ij}$ be the extension of $\widehat{A}^{ij}$ which follows the same rule as $a^{ij}$.

By a direct computation, $(\tilde{u},\tilde{p})$ satisfies
\[
\partial_t \tilde{u}-\overline{a}^{ij}D_{ij} \tilde{u}-\overline{b}^i D_i \tilde{u}-\overline{c}\tilde{u}+\nabla \tilde{p}=\tilde{f},\quad \Div (\overline{\widehat{A}}\tilde{u})=\tilde{g}\quad \text{in } Q_{R},\]
and
\[    D_d \tilde{u}^k=\tilde{u}^d=0\quad \text{on } (-R^2,0]\times B_{R}' \times \{0\}.\]

Choose a mollifier $\varphi_\varepsilon$ which is even with respect to the $x_d$ variable. 
Define
\[
     \tilde{u}^{(\varepsilon)}(t,x',x_d)=\int_{Q_\varepsilon} \varphi_\varepsilon(s,y',y_d)\tilde{u}(t-s,x-y',x_d-y_d)\myd{y}{ds}
\]
for $(t,x) \in (-(R^*)^2+\varepsilon^2,0)\times B_{R^*-\varepsilon}$, $\varepsilon\in(0,1)$.

By a change of variables, {we have $(\tilde{u}^d)^{(\varepsilon)}(t,x',0)=0$} since $\varphi_\varepsilon$ is a symmetric mollifier with respect to the $x_d$ variable. Similarly, we have $D_d (\tilde{u}^{(\varepsilon)})^k(t,x',0)=0$ for $k=1,2,\dots,d-1$. Hence $\tilde{u}^{(\varepsilon)}$ satisfies the Lions boundary conditions on $(-R^2+\varepsilon^2,0]\times B_{R}'\times \{0\}$.

Choose $\beta,\gamma>0$ and $R_0>0$ so that Lemma \ref{lem:regularity-lemma-Stokes}, Theorems \ref{thm:solvability-modified-divergence} and \ref{thm:regularity-a-priori-local} hold. We can further extend $\widehat{A}$ to the whole domain $\mathbb{R}^d$ so that it satisfies Assumption \ref{assump:delta-N} $(\beta)$. Indeed, choose $\eta \in C_0^\infty (B_{R_0})$ so that $\eta =1$ in $B_{R_0/2}$ and $0\leq \eta \leq 1$. Then if we define
\[ \widehat{A}_\eta = \eta \widehat{A} + (1-\eta)I, \]
then $\widehat{A}_\eta = \widehat{A}$ in $B_{R_0/2}$ and 
\[ \norm{\widehat{A}_\eta -I}{\Leb{\infty}}\leq \norm{\eta (\widehat{A}-I)}{\Leb{\infty}(B_{R_0})}\leq \beta. \]

For simplicity, we drop a tilde and overbar in the notation. If we write $\theta=\sqrt{\frac{2}{3}}$, then note that 
\[ R/2= 3\theta^2 R/4<\theta^2 R<\theta R<R.\]
  Choose a cut-off function $\psi \in C_0^\infty (Q_{R})$ so that $\psi=1$ in $Q_{\theta^2 R}$ and $\psi=0$ outside $Q_{\theta R}$. By taking mollification in $(t,x)$, for $\varepsilon$ small, we have
\[  \partial_t u^{(\varepsilon)}-\mathcal{L}u^{(\varepsilon)}+\nabla p^{(\varepsilon)}=f^{(\varepsilon)}+{f}^\varepsilon_e,\quad \Div(\widehat{A}u^{(\varepsilon)})=g^{(\varepsilon)}+g^\varepsilon_e\quad
\text{ in } Q_{{\theta^2 R}}. \]
where
\begin{align*}
\mathcal{L}u&=a^{ij}(t,x)D_{ij} u +b^i (t,x) D_i u + cu,\\
{f}^\varepsilon_e&=\left[(\mathcal{L}u)^{(\varepsilon)}-\mathcal{L}u^{(\varepsilon)}\right]1_{Q_{\theta^2 R}},\quad
{g}^\varepsilon_e=\Div[(\widehat{A}u^{(\varepsilon)}-(\widehat{A}u)^{(\varepsilon)})\psi].
\end{align*}
Since $g_e^\varepsilon$ satisfies the compatibility condition in Theorem \ref{thm:solvability-modified-divergence}, it follows from Theorem \ref{thm:solvability-modified-divergence} that there exists $(u_1^\varepsilon,p_1^\varepsilon)\in \oSob{1,2}{q_0}((-R^2,0)\times\mathbb{R}^d_+)^d\times \Sob{0,1}{q_0}((-R^2,0)\times\mathbb{R}^d_+)$ satisfying 
\begin{equation}\label{eq:perturbed-Stokes-II}
\left\{
\begin{alignedat}{2}
\partial_t u -a^{ij}(t,x)D_{ij} u +b^i(t,x)D_i u + c(t,x)u+\nabla p &=f_e^\varepsilon &&\quad \text{in } (-R^2,0)\times\mathbb{R}^d_+,\\
\Div (\widehat{A} u)&=g_e^\varepsilon&&\quad \text{in } (-R^2,0)\times\mathbb{R}^d_+\\ 
\end{alignedat}
\right.
\end{equation}
under the Lions boundary conditions \eqref{eq:Lions}. Moreover, we have
\begin{equation}\label{eq:u-1-epsilon}
   \norm{u_1^\varepsilon}{\Sob{1,2}{q_0}((-R^2,0)\times\mathbb{R}^d_+)}\leq N \left(\norm{f_e^\varepsilon}{\Leb{q_0}(Q_{\theta^2 R}^+)}+\norm{[\widehat{A}u^{(\varepsilon)}-(\widehat{A}u)^{(\varepsilon)}]\psi}{\Sob{1,2}{q_0}(Q_{\theta R}^+)}\right).
\end{equation}
By Lemma \ref{lem:regularity-lemma-Stokes}, $u_1^\varepsilon \in \Sob{1,2}{s,q}((-R^2,0)\times\mathbb{R}^d_+)$. Define $u_2^\varepsilon=u^{(\varepsilon)}-u_1^\varepsilon$ and $p_2^\varepsilon=p^{(\varepsilon)}-p_1^\varepsilon$. Then by Theorem \ref{thm:regularity-a-priori-local}, we have
\begin{equation}
            \label{eq:approximate}
\begin{aligned}
    & \norm{D^2 u_2^\varepsilon}{\Leb{s,q}(Q_{R/2}^+)}\leq N\left( \norm{u_2^\varepsilon}{\Leb{s,1}(Q_{\theta^2 R}^+)}+\norm{f^{(\varepsilon)}}{\Leb{s,q}(Q_{\theta^2 R}^+)}+\norm{Dg^{(\varepsilon)}}{\Leb{s,q}(Q_{\theta^2 R}^+)}\right)\\
 &\leq N\left( \norm{u^{(\varepsilon)}}{\Leb{s,1}(Q_{\theta^2 R}^+)}+\norm{u_1^\varepsilon}{\Leb{s,1}(Q_{\theta^2 R}^+)}+\norm{f^{(\varepsilon)}}{\Leb{s,q}(Q_{\theta^2 R}^+)}+\norm{Dg^{(\varepsilon)}}{\Leb{s,q}(Q_{\theta^2 R}^+)}\right),
\end{aligned}        
\end{equation}
where $N=N(d,s,q,\nu,R,R_0)>0$. By the Sobolev embedding theorem and \eqref{eq:u-1-epsilon}, we have 
\begin{equation}\label{eq:u-1-epsilon-perturb}
   u\in \Leb{s,1}(Q_R^+)^d\quad\text{and}\quad  \norm{u_1^\varepsilon}{\Leb{s,1}(Q_R^+)}\rightarrow 0 
\end{equation}
as $\varepsilon\rightarrow 0$. Since $f \in \Leb{s,q}(Q_{R}^+)^d$ and $g \in \Sob{0,1}{s,q}(Q_{R}^+)$, it follows from \eqref{eq:approximate} that there exists a constant $N>0$ independent of $\varepsilon$ such that 
\[  \sup_{\varepsilon>0} \norm{D_{ij} u_2^\varepsilon}{\Leb{s,q}(Q_{R/2}^+)}\leq N.\]
Hence by the weak compactness in $\Leb{s,q}(Q_{R/2}^+)$, there exists a subsequence $\{D_{ij} u_2^{\varepsilon_k}\}$ of $\{D_{ij} u_2^\varepsilon\}$ which converges weakly to a function $v_{ij}$ in $\Leb{s,q}(Q_r^+)$. 

On the other hand, we have $D^2 u^{(\varepsilon)}\rightarrow D^2 u$ strongly in $\Leb{q_0}(Q_{\theta R}^+)$ and $D^2 u_1^\varepsilon\rightarrow 0$ strongly in $\Leb{q_0}(Q_{\theta R}^+)$ as $\varepsilon\rightarrow 0+$. Hence it follows that $D^2 u_2^\varepsilon \rightarrow D^2 u$ strongly in $\Leb{q_0}(Q_{\theta R}^+)$, which implies that $D^2 u=v$ in $Q_r^+$. Therefore we get the desired inequality by taking liminf in \eqref{eq:approximate}. This completes the proof of Theorem \ref{thm:main-theorem-perturb}.
\end{proof}

\section{Proof of Theorem \ref{thm:A}} \label{sec:5}
This section is devoted to a proof of Theorem \ref{thm:A}. We first show that for each $x_0 \in \partial\Omega$, we can transform the equation into a perturbed Stokes equation on the half-space so that we can apply Theorem \ref{thm:main-theorem-perturb}. Then the desired result follows by a covering argument. 

\begin{proof}[Proof of Theorem \ref{thm:A}]
Fix $x_0 \in \partial\Omega$ and take $z_0=(0,x_0)$. Define a flattening map $\Phi:\mathbb{R}^d\rightarrow\mathbb{R}^d$ by \eqref{eq:transformation-1} which becomes a local diffeomorphism with inverse $\Psi: B_{R_1^*}(x_0)\rightarrow\mathbb{R}^d$ satisfying $\Psi(B_{R_1^*}(x_0)\cap\Omega)\subset \mathbb{R}^d_+$ for some $R_1^*\in (0,R_0^*)$. There exists $R_2^*>0$ depending on $R_0^*$ and $M$ such that 
\[  \Phi(B_{R_2^*})\subset B_{R_1^*}(x_0)\quad \text{and}\quad \Phi(B_{R_2^*}^+)\subset \widehat{B}_{R_1^*}(x_0),\]
where we assumed, without loss of generality, $0=y_0=\Psi(x_0)$. Define 
\[  \tilde{u}(t,y)=\left(\frac{\partial x}{\partial y}\right)^T \hat{u}(t,y),\quad \text{where } \hat{u}(t,y)=u(t,\Phi(y))\]
for $(t,y) \in Q_{R_2^*}^+$. Write
\begin{equation}\label{eq:delta-E}
\frac{\partial x_i}{\partial y_j}(y)=\delta_{ij} + E_{ij}(y),
\end{equation}
where $E_{ij}(y)$ satisfies
\begin{equation}\label{eq:E-estimate}
 |E_{ij}(y)|\leq N\theta\quad \text{and}\quad |D^k E_{ij}(y)|\leq N,\quad k=1,2 
\end{equation}
for some constant $N=N(d,M)>0$ and for all $y\in B_{R_0^*/2}$ by Lemma \ref{lem:estimates-boundary}.	

Since 
\[ \hat{u}(t,y)=\left(\frac{\partial y}{\partial x}\right)^T \tilde{u}(t,y),\]
differentiating in $y_k$, we get 
\[
   (D_ku)(t,\Phi(y))=-(D_lu)(t,\Phi(y))E_{lk}(y)+\left(\frac{\partial y}{\partial x}\right)^T \frac{\partial \tilde{u}}{\partial y_k} + \frac{\partial}{\partial y_k}\left(\frac{\partial y}{\partial x}\right)^T\tilde{u}(t,y)\]
and it follows from \eqref{eq:delta-E} and \eqref{eq:E-estimate} that 
\begin{equation}\label{eq:tilde-Du}
\norm{(Du)(\cdot,\Phi(\cdot))}{\Leb{s,q}(Q_r^+)}
\leq N\theta \norm{(Du)(\cdot,\Phi(\cdot))}{\Leb{s,q}(Q_r^+)}+N\norm{D\tilde{u}}{\Leb{s,q}(Q_r^+)} + N\norm{\tilde{u}}{\Leb{s,q}(Q_r^+)} 
\end{equation}
for some constant $N=N(d,s,q,M)>0$. A similar computation gives
\begin{equation}\label{eq:tilde-D2u}
\begin{aligned}
\norm{(D^2 u)(\cdot,\Phi(\cdot))}{\Leb{s,q}(Q_r^+)}
&\leq N\theta \norm{(D^2u)(\cdot,\Phi(\cdot))}{\Leb{s,q}(Q_r^+)}+N\norm{(Du)(\cdot,\Phi(\cdot))}{\Leb{s,q}(Q_r^+)} \\
&\relphantom{=}+ N\left(\norm{D^2\tilde{u}}{\Leb{s,q}(Q_r^+)} +\norm{D\tilde{u}}{\Leb{s,q}(Q_r^+)} +\norm{\tilde{u}}{\Leb{s,q}(Q_r^+)}\right)
\end{aligned}
\end{equation}
for some constant $N=N(d,s,q,M)>0$. By choosing $\theta$ sufficiently small in \eqref{eq:tilde-Du} and \eqref{eq:tilde-D2u}, we have
\begin{align*}
\norm{(D^2u)(\cdot,\Phi(\cdot))}{\Leb{s,q}(Q_r^+)}&\leq N\left(\norm{D^2\tilde{u}}{\Leb{s,q}(Q_r^+)} +\norm{D\tilde{u}}{\Leb{s,q}(Q_r^+)} +\norm{\tilde{u}}{\Leb{s,q}(Q_r^+)}\right)
\end{align*}
for some constant $N=N(d,s,q,M)>0$. Then by the Poincar\'e inequality and the interpolation inequality \eqref{eq:Du-interpolation-L1}, we get 
\begin{equation}\label{eq:Hessian-u-sol-change}
\begin{aligned}
&\norm{(D^2u)(\cdot,\Phi(\cdot))}{\Leb{s,q}(Q_r^+)}\\
&\leq N\left(\norm{D^2\tilde{u}}{\Leb{s,q}(Q_r^+)} +(r+1)\norm{D\tilde{u}}{\Leb{s,q}(Q_r^+)} +r^{d(1/q-1)}\norm{\tilde{u}}{\Leb{s,1}(Q_r^+)}\right)\\
&\leq N\norm{D^2\tilde{u}}{\Leb{s,q}(Q_r^+)} + N\left(r^{-1+d(1/q-1)} +(r+1)^{\theta/(1-\theta)}\right)\norm{u(\cdot,\Phi(\cdot))}{\Leb{s,1}(Q_r^+)}
\end{aligned}
\end{equation}
for some constant $N=N(d,s,q,M)>0$. 

By \eqref{eq:no-penetration} and Proposition \ref{prop:Lions-flatterning}, one can check that if $(u,p)$ satisfies \eqref{eq:Stokes} in $Q_{R_0^*}(0,x_0) \cap (\mathbb{R}\times \Omega)$ with the boundary condition  \eqref{eq:general-Lions-boundary} on $(-(R_0^*)^2,0] \times (B_{R_0^*}(x_0) \cap \partial\Omega)$, then we have
\begin{equation*}
\partial_t \hat{u}-\tilde{a}^{ij}D_{ij}\hat{u}-\tilde{b}^i D_i\hat{u}+\left(\frac{\partial y}{\partial x}\right)^T\nabla_y \hat{p}=\hat{f}\quad \text{in } Q_{R_2^*}^+,
\end{equation*}
where
\begin{align*}
\tilde{a}^{ij}(t,y)=a^{lm}(t,\Phi(y))\frac{\partial y_i}{\partial x_l}\frac{\partial y_j}{\partial x_m}(\Phi(y)),\quad
\tilde{b}^{i}(t,y)=a^{lm}(t,\Phi(y))\frac{\partial^2 y_i}{\partial x_l \partial x_m}(\Phi(y)).
\end{align*}
Note that 
\[  \partial_t \tilde{u}-\mathcal{L}\tilde{u}+\left[\mathcal{L}\tilde{u}-\left(\frac{\partial x}{\partial y}\right)^T \mathcal{L}\hat{u}\right]+\nabla_y \hat{p}=\left(\frac{\partial x}{\partial y}\right)^T \hat{f}.
\]
A direct computation and Lemma \ref{lem:estimates-boundary} give
\begin{align*}
\partial_t \tilde{u}-\tilde{a}^{ij}D_{ij}\tilde{u}-\mathsf{b}^iD_i \tilde{u}-\mathsf{c} \tilde{u}+\nabla \hat{p}=\left(\frac{\partial x}{\partial y}\right)^T \hat{f}\quad \text{in } Q_{R_2^*}^+,
\end{align*}
where 
\[  |\mathsf{b}(t,y)|+|\mathsf{c}(t,y)|\leq N.\]
Moreover, since $\Div u=g$ in $Q_{R_0^*}(0,x_0)\cap (\mathbb{R}\times \Omega)$, one can easy to see that 
\begin{align*}
\Div_y\left(\widehat{C}\tilde{u}\right)&=\hat{g}\left|\det \left(\frac{\partial x}{\partial y}\right)\right|\quad \text{in } Q_{R_2^*}^+,
\end{align*}
where 
\begin{equation}\label{eq:hat-C-identity}
    \widehat{C}(y)= \left|\det \left(\frac{\partial x}{\partial y}\right)\right|\left(\frac{\partial y}{\partial x}\right)\left(\frac{\partial y}{\partial x}\right)^T(\Phi(y))
\end{equation}
is a symmetric matrix. Then it is easily seen from \eqref{eq:delta-E}, \eqref{eq:E-estimate}, and \eqref{eq:hat-C-identity} that 
\[ \norm{\hat{C}-I}{\Leb{\infty}(B_{R_2^*}^+)}\leq N\theta\quad\text{and}\quad \sum_{k=1}^2 \norm{D^k\hat{C}}{\Leb{\infty}(B_{R_2^*}^+)}\leq N,\]
where $N=N(d,M)>0$. 

Now we need to check the conditions of Theorem \ref{thm:main-theorem-perturb}. By Lemma \ref{lem:estimates-boundary}, one can check that 
\begin{equation}\label{eq:coefficient-condition-check}
 \sup_{i,j} \fint_{Q_r(x_1)} |\tilde{a}^{ij}(t,x)-(\tilde{a}^{ij})_{B_r(x_1)}(t)| \myd{x}dt \leq N_0(\theta +\gamma) 
\end{equation}
for any $x_1 \in B_{R_2^*/2}$ and $0<r\leq R_2^*/2$, where $N_0$ is independent of $\theta$ and $\gamma$. We may change the values of $\tilde{a}^{ij}$ outside $Q_{R_2^*/2}$ and extend them to $\mathbb{R}^{d+1}$ so that $\tilde{a}^{ij}$ satisfies boundedness, uniform ellipticity conditions, and \eqref{eq:coefficient-condition-check} with possibly different $N_0$ and small radius. Indeed, this can be shown by using
\[ \tilde{a}^{ij}\eta(\cdot/R_2^*)+\delta^{ij}(1-\eta(\cdot/R_2^*)),\]
where $\eta \in C_0^\infty(B_{1/4})$ satisfying $\eta=1$ in $B_{1/8}$. 

By Theorem \ref{thm:main-theorem-perturb}, for $0<r_1\leq \frac{1}{4}\min\{1,R_0^*,R_2^*\}$, we have 
\begin{equation}\label{eq:D2-tilde}
\norm{D^2\tilde{u}}{\Leb{s,q}(Q_{r_1/2}^+)}\leq N\left(\norm{\tilde{u}}{\Leb{s,1}(Q_{r_1}^+)}+\norm{\hat{f}}{\Leb{s,q}(Q_{r_1}^+)}+\norm{\hat{g}}{\Sob{0,1}{s,q}(Q_{r_1}^+)} \right)
\end{equation}
for some constant $N=N(d,s,q,\nu,r_1,R_0)>0$. By \eqref{eq:Hessian-u-sol-change} and \eqref{eq:D2-tilde}, we have
\begin{align*}
\norm{(D^2u)(\cdot,\Phi(\cdot))}{\Leb{s,q}(Q_{r_1/2}^+)}&\leq N\left(\norm{\hat{u}}{\Leb{s,1}(Q_{r_1}^+)}+\norm{\hat{f}}{\Leb{s,q}(Q_{r_1}^+)}+\norm{\hat{g}}{\Sob{0,1}{s,q}(Q_{r_1}^+)}\right)
\end{align*}
for some constant $N=N(d,s,q,\nu,M,r_1,R_0)>0$. Hence a change of variables gives 
\begin{equation}
    \label{eq:D2u-estimates-small-scale}
\begin{aligned}
    &\norm{D^2u}{\Leb{s,q}((-r_1^2/4,0)\times \Phi(B_{r_1/2}^+))}\\
&\leq N\left(\norm{u}{\Leb{s,1}((-r_1^2,0)\times \Phi(B_{r_1}^+))}+\norm{f}{\Leb{s,q}((-r_1^2,0)\times \Phi(B_{r_1}^+))}+\norm{g}{\Sob{0,1}{s,q}((-r_1^2,0)\times \Phi(B_{r_1}^+))}\right)
\end{aligned}
\end{equation}
for some constant $N=N(d,s,q,\nu,M,r_1,R_0)>0$. 

\begin{figure}[h]
\includegraphics[width=0.9\textwidth]{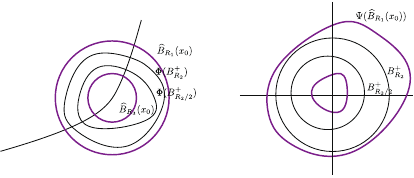}
\caption{}\label{fig:1}
\end{figure}

Now for $R_1>0$ satisfying $0<R_1\leq \frac{1}{4} \min \{1,R_0^*,R_2^*\}$, choose $R_2>0$ so that $R_2\in(0,R_1)$ and $B_{R_2}^+\subset \Psi(\widehat{B}_{R_1}(x_0))$. Next, choose $R_3 \in (0,R_2/2)$ so that $\widehat{B}_{R_3}(x_0)\subset \Phi(B_{R_2/2}^+)$ (see Figure \ref{fig:1}). Then by \eqref{eq:D2u-estimates-small-scale}, we have 
\begin{align*}
&\norm{D^2u}{\Leb{s,q}(\widehat{Q}_{R_3}(z_0))}\leq \norm{D^2u}{\Leb{s,q}((-R_2^2/4,0)\times \Phi(B_{R_2/2}^+))}\\
&\leq N\left(\norm{u}{\Leb{s,1}((-R_2^2,0)\times \Phi(B_{R_2}^+))}+\norm{f}{\Leb{s,q}((-R_2^2,0)\times \Phi(B_{R_2}^+))}+\norm{g}{\Sob{0,1}{s,q}((-R_2^2,0)\times \Phi(B_{R_2}^+))}\right)\\
&\leq N\left(\norm{u}{\Leb{s,1}(\widehat{Q}_{R_1}(z_0))}+\norm{f}{\Leb{s,q}(\widehat{Q}_{R_1}(z_0))}+\norm{g}{\Sob{0,1}{s,q}(\widehat{Q}_{R_1}(z_0))}\right)
\end{align*}
for some constant $N=N(d,s,q,\nu,M,R_1,R_0)>0$. By the translation invariance in time, we can shift the temporal center of the cylinder. Since the constant $N$ does not depend on the point $z_0$, the desired result follows by a covering argument. This completes the proof of Theorem \ref{thm:A}.
\end{proof}

\section{Proof of Theorem \ref{thm:B}} \label{sec:6}

This section is devoted to the proof of Theorem \ref{thm:B}. For simplicity, we introduce additional notation. We say that $A\apprle B$ if there exists a constant $N>0$ independent of $A$ and $B$ such that $A\leq NB$.

For $x_0 \in \partial\Omega$, there exist $r=r_{x_0}>0$ and $\phi$ such that 
\[  B_r(x_0)\cap \Omega = \{ x\in B_r(x_0): x_d>\phi(x')\}.\]

Define $\Psi$ by 
\[ y_k=\Psi^k(x)=x_k,\quad k=1,\dots,d-1,\quad \text{and} \quad y_d =\Psi^d(x)=x_d - \phi(x')\]
and $\Phi$ by
\[ x_k=\Phi^k(y)=y_k,\quad k=1,\dots,d-1,\quad \text{and}\quad x_d=\Phi^d(y)=y_d+\phi(y').\]
Then $\Phi=\Psi^{-1}$,
\[  D\Phi(y)= \begin{bmatrix}
1 & 0 & \cdots & 0 \\
0 & 1 & \cdots & 0 \\
\vdots & \vdots & & \vdots \\
D_1\phi(y') & D_2 \phi(y') & \cdots & 1
\end{bmatrix}\,\,
\text{and}
\,\,
D\Psi(x)=\begin{bmatrix}
1 & 0 & \cdots & 0 \\
0 & 1 & \cdots & 0 \\
\vdots & \vdots & & \vdots \\
-D_1\phi(x') & -D_2 \phi(x') & \cdots & 1
\end{bmatrix}.
\]
Moreover, we have $\det D\Phi=\det D\Psi =1$. Define 
\[
\mathcal{L}\hat{u}=\hat{a}^{ij}D_{ij} \hat{u}+\hat{b}^i D_i \hat{u},
\]
where 
\begin{align*}
  \hat{a}^{ij}(t,y)=a^{kl}(t,\Phi(y))\Psi^i_{x_l}(\Phi(y))\Psi^j_{x_k}(\Phi(y)), \quad
  \hat{b}^i(t,y)=a^{kl}(t,\Phi(y))\frac{\partial^2 \Psi^i}{\partial x_k \partial x_l}(\Phi(y)),
\end{align*}
and
\[ \begin{gathered}
\hat{u}(t,y)=u(t,\Phi(y)),\quad \hat{p}(t,y)=p(t,\Phi(y)),\\
\hat{f}(t,y)=f(t,\Phi(y)),\quad \text{and}\quad \hat{g}(t,y)=g(t,\Phi(y)).
\end{gathered}  \]

Then it is easily seen that $(\hat{u},\hat{p})$ satisfies
\begin{equation}\label{eq:usual-flattening}
 \partial_t \hat{u} -\mathcal{L}\hat{u}+\nabla \hat{p}=\hat{f}+(D_d \hat{p}) (D_1 \phi, D_2\phi,\dots,D_{d-1}\phi,0)^T \quad\text{and}\quad \Div(\widehat{C} \hat{u})=\hat{g},
\end{equation}
where 
\[ \widehat{C}(y)=\begin{bmatrix}
1& 0 & \cdots & 0 \\
0 & 1 & \cdots & 0 \\
\vdots & \vdots & & \vdots \\
-D_1\phi(y') & -D_2\phi(y') & \cdots & 1
\end{bmatrix}.
\]
Moreover, $\hat{g}$ satisfies $\partial_t \hat{g}=\Div \widehat{G}$ in $Q_R^+$, where 
\begin{equation}\label{eq:G-hat-defn}
\widehat{G}(t,y)=G(t,\Phi(y))[D\Phi(y)]^{-T}.
\end{equation}

To interpret the boundary condition \eqref{eq:Navier-curved}, we recall that $(2\mathbb{D}u)_{ij}=D_ju^i +D_i u^j$ and we choose $\tau^k=e_k + (D_k \phi)e_d$. The boundary condition \eqref{eq:Navier-curved} implies that 
\begin{equation}\label{eq:ud-expression}
u^d(t,y',\phi(y'))=\sum_{l=1}^{d-1} u^l(t,y',\phi(y'))(D_l \phi)(y')
\end{equation}
and 
\begin{equation*}\label{eq:Navier-translation}
\begin{aligned}
0&=\sum_{i,j=1}^d\left[(2\mathbb{D}u)_{ij}n_j \tau_i^k + (\mathcal{A} u)^i \tau_i^k\right]\\
&=\sum_{j=1}^d\left[(2\mathbb{D}u)_{kj}n_j + (2\mathbb{D}u)_{dj} n_j (D_k\phi)\right]+(\mathcal{A}u)^k + (\mathcal{A}u)^d D_k \phi.
\end{aligned}
\end{equation*}
for all $k=1,\dots,d-1$. Then by the representation of $n$, we have
\begin{equation}\label{eq:Navier-translation-2}
\begin{aligned}
0&=\sum_{j=1}^{d-1}\left[(2\mathbb{D}u)_{kj} (D_j \phi)+(2\mathbb{D}u)_{dj}(D_j \phi)(D_k\phi)\right]\\
&\relphantom{=}+(2\mathbb{D}u)_{kd}(-1) + (2\mathbb{D}u)_{dd}(-1)(D_k \phi)+\sqrt{1+|D\phi|^2} ((\mathcal{A}u)^k +(\mathcal{A}u)^d D_k \phi).
\end{aligned}
\end{equation}

On the other hand, by differentiating \eqref{eq:ud-expression} with respect to $y_k$, we get
\begin{equation}\label{eq:ud-k}
\begin{aligned}
D_k u^d =-(D_d u^d)(D_k \phi)+\sum_{l=1}^{d-1} \left[(D_k u^l)(D_l \phi)+(D_d u^l)(D_k \phi)(D_l \phi)+u^l (D_{kl}\phi)\right]
\end{aligned}
\end{equation}
evaluated at $(t,y',\phi(y'))$ for $k=1,\dots,d-1$. Taking the difference of \eqref{eq:Navier-translation-2} and \eqref{eq:ud-k}, we have for $k=1,\dots,d-1$,
\begin{align*}
D_d u^k&=-(D_d u^d)(D_k \phi)-\sum_{l=1}^{d-1}\left[ (D_k u^l)(D_l \phi)+(D_d u^l)(D_k\phi)(D_l\phi)+u^l (D_{kl}\phi)\right]\\
&\relphantom{=}+\sum_{l=1}^{d-1}\left[(2\mathbb{D}u)_{kl} (D_l \phi)+(2\mathbb{D}u)_{dl}(D_l \phi)(D_k\phi)\right]+\sqrt{1+|D\phi|^2} ((\mathcal{A}u)^k +(\mathcal{A}u)^d D_k \phi)\\
&=-(D_d u^d)(D_k \phi)+\sum_{l=1}^{d-1} [(D_l u^k)(D_l \phi)+(D_l u^d)(D_k \phi)(D_l \phi)-u^l(D_{kl}\phi) ]\\
&\relphantom{=}+\sqrt{1+|D\phi|^2} ((\mathcal{A}u)^k +(\mathcal{A}u)^d D_k \phi).
\end{align*}

Define
\begin{align*}
h_0(t,y')&= \sum_{l=1}^{d-1} u^l(t,y',\phi(y')) (D_l\phi)(y'),\\
h_k(t,y')&=-(D_d u^d)(t,y',\phi(y'))(D_k \phi)(y')-\sum_{l=1}^{d-1}u^l(t,y',\phi(y'))(D_{kl}\phi)(y')\\
&\relphantom{=}+\sqrt{1+|D\phi(y')|^2} \left((\mathcal{A}(t,y',\phi(y'))u)^k(t,y',\phi(y')) +(\mathcal{A}(t,y',\phi(y'))u)^d(t,y',\phi(y') D_k \phi(y')\right)\\
&\relphantom{=}+\sum_{l=1}^{d-1} [(D_l u^k)(t,y',\phi(y'))(D_l \phi)(y')+(D_l u^d)(t,y',\phi(y'))(D_k \phi)(y')(D_l \phi)(y')]
\end{align*}
for $k=1,\dots,d-1$.

Then $u$ satisfies 
\[  u^d(t,y',\phi(y'))=h_0(t,y')\quad \text{and}\quad D_d u^k(t,y',\phi(y'))=h_k(t,y'),\quad k=1,\dots,d-1,\]
and so $\hat{u}$ satisfies
\begin{equation}\label{eq:transformed-A}
\left\{
\begin{aligned}
&\partial_t \hat{u}-(\tilde{a}^{ij}D_{ij}\hat{u}+\tilde{b}^i D_i \hat{u})+\nabla \hat{p}=\hat{f}+(D_d \hat{p})(D_1\phi,D_2\phi,\dots,D_{d-1}\phi,0)^T,\\
&\Div(\widehat{C}\hat{u})=\hat{g},\\
\end{aligned}
\right.
\end{equation}
in $Q_{R}^+$ and 
\begin{equation}\label{eq:transformed-B}
\hat{u}^d = h_0,\quad D_d\hat{u}^k=h_k
\end{equation}
on $(-R^2,0]\times B_{R}'\times \{0\}$.

We may assume that $R\leq 1$. For $0<r<R$, choose $\zeta \in C_0^\infty (Q_R)$ so that $\zeta=1$ in $Q_r$ and 
\[ |\partial_t \zeta |+|D^2 \zeta|\apprle \frac{1}{(R-r)^2}\quad \text{and}\quad |D\zeta|\apprle \frac{1}{R-r}. \]

Define 
\[  \tilde{u}=\hat{u}\zeta\quad\text{and}\quad \tilde{p}=\hat{p}\zeta. \]
Then $(\tilde{u},\tilde{p}) \in \Sob{1,2}{s,q}((-R^2,0)\times \mathbb{R}^d_+)^d\times \Sob{0,1}{s,q}((-R^2,0)\times \mathbb{R}^d_+)$ and satisfy
\begin{equation*}
\left\{\begin{aligned}
\partial_t \tilde{u}-\mathcal{L}\tilde{u}+\nabla \tilde{p}&=\tilde{f},\\
\Div(\widehat{C}\tilde{u})&=\tilde{g}
\end{aligned}
\right.
\end{equation*}
in $(-R^2,0)\times \mathbb{R}^d_+$ subject to the boundary condition
\[
   \tilde{u}^d=\zeta h_0,\quad D_d \tilde{u}^k=\zeta h_k + \hat{u}^k D_d \zeta,\quad k=1,\dots,d-1\quad\text{on }(-R^2,0]\times \mathbb{R}^{d-1}\times \{0\}.
\]
Here 
\begin{equation}\label{eq:f-tilde-g-tilde}
\begin{aligned}
\tilde{f}&=\zeta(\hat{f}+(D_d \hat{p})(D_1 \phi,D_2\phi,\dots, D_{d-1}\phi,0)^T)+2(D_i\zeta) (D_{j}\hat{u}) a^{ij} + \hat{p}\nabla \zeta\\
&\relphantom{=}+\hat{u}(\partial_t \zeta -\mathcal{L}\zeta), \\
\tilde{g}&= \hat{g} \zeta + (\widehat{C}\hat{u})\cdot \nabla \zeta.
\end{aligned}
\end{equation}
Note that $\tilde{u}(-R^2,\cdot)=0$ by the choice of $\zeta$.

For $\lambda> 0$, let $v=(v_1,\dots,v_d) \in \Sob{1,2}{s,q}((-R^2,\infty)\times \mathbb{R}^d_+)^d$ be the unique solution to the following Cauchy problems:
\begin{equation*}
\left\{
\begin{aligned}
\partial_t v_k-\Delta v_k+\lambda v_k&=0 \\
D_d v_k &=\zeta h_k +\hat{u}^k D_d \zeta
\end{aligned}
\right.
\quad\text{and}\quad \left\{
\begin{aligned}
\partial_t v_d-\Delta v_d+\lambda v_d&=0 \\
v_d &=\zeta h_0
\end{aligned}
\right.
\end{equation*}
for $k=1,\dots,d-1$ and $v(-R^2,\cdot)=0$. 

Then for $\varepsilon\in(0,1)$ and $\lambda \geq 1$, by using interpolation inequalities and the trace theorem, we have
\begin{equation}\label{eq:v-estimate}
\begin{aligned}
&\lambda \norm{v}{\Leb{s,q}((-R^2,\infty)\times \mathbb{R}^d_+)}+\sqrt{\lambda}\norm{Dv}{\Leb{s,q}((-R^2,\infty)\times \mathbb{R}^d_+)}\\
&\relphantom{=}+\norm{D^2 v}{\Leb{s,q}((-R^2,\infty)\times \mathbb{R}^d_+)}+\norm{\partial_t v}{\Leb{s,q}((-R^2,\infty)\times \mathbb{R}^d_+)}\\
&\leq N(\varepsilon + \norm{D\phi}{\Leb{\infty}(B_R')})\left(\norm{\partial_t \hat{u}}{\Leb{s,q}(Q_R^+)}+\norm{D^2\hat{u}}{\Leb{s,q}(Q_R^+)} \right)\\
&\relphantom{=}+\frac{N\lambda^b}{\varepsilon^{2/(q-1)+1}(R-r)^2}\norm{\hat{u}}{\Leb{s,q}(Q_R^+)},
\end{aligned}
\end{equation}
where $N$ is independent of $\varepsilon$, $\lambda$, $r$, $R$, and  $b=2-1/q$. 
In particular, if $\varepsilon,\theta,R \in (0,1)$, $\sqrt{\lambda}=(R-r)^{-1}\theta^{-1}$, then \eqref{eq:v-estimate} can be simplified as below:
\begin{equation}\label{eq:v-estimate-simplified}
\begin{aligned}
&\lambda \norm{v}{\Leb{s,q}((-R^2,\infty)\times \mathbb{R}^d_+)}+\sqrt{\lambda}\norm{Dv}{\Leb{s,q}((-R^2,\infty)\times \mathbb{R}^d_+)}\\
&\relphantom{=}+\norm{D^2 v}{\Leb{s,q}((-R^2,\infty)\times \mathbb{R}^d_+)}+\norm{\partial_t v}{\Leb{s,q}((-R^2,\infty)\times \mathbb{R}^d_+)}\\
&\leq N(\varepsilon + \norm{D\phi}{\Leb{\infty}(B_R')})\left(\norm{\partial_t \hat{u}}{\Leb{s,q}(Q_R^+)}+\norm{D^2\hat{u}}{\Leb{s,q}(Q_R^+)}\right)\\
&\relphantom{=}+\frac{N_{\theta,\varepsilon}}{(R-r)^{2+2b}} \norm{\hat{u}}{\Leb{s,q}(Q_R^+)},
\end{aligned}
\end{equation}
where $N_{\theta,\varepsilon}$ depends on $\theta$ and $\varepsilon$.

To proceed further, we prove the following lemma.
\begin{lemma}\label{lem:lambda-estimate}
Let $\lambda>0$ and $s,q\in(1,\infty)$. Under Assumption \ref{assump:VMO} $(\gamma)$, if $(u,p) \in \oSob{1,2}{s,q}(\Omega_T)^d\times \Sob{0,1}{s,q}(\Omega_T)$ satisfies 
\begin{equation*}
\left\{
\begin{aligned}
\partial_t u -(a^{ij}D_{ij}u+b^i D_i u) +\nabla p +\lambda u&=f \\
\Div u &=g
\end{aligned}
\right.
\end{equation*}
in $\Omega_T$, $u(0,\cdot)=0$, and the Lions boundary conditions
\[ u^d =0,\quad D_d u^k=0,\quad k=1,\dots,d-1\quad \text{on } (0,T)\times \mathbb{R}^{d-1}\times \{0\}\]
for some $f\in \Leb{s,q}(\Omega_T)^d$ and $g\in \mathring{\mathcal{H}}^1_{s,q}(\Omega_T)$ satisfying $g_t = \Div G+ h$ for some $G,h\in \Leb{s,q}(\Omega_T)$, then there exists a constant $N=N(d,s,q,\nu,R_0)>0$ such that 
\begin{align*}
&\lambda \norm{u}{\Leb{s,q}((0,T)\times \mathbb{R}^d_+)}+\sqrt{\lambda}\norm{Du}{\Leb{s,q}((0,T)\times \mathbb{R}^d_+)}\\
&\relphantom{=}+\norm{\partial_t u}{\Leb{s,q}((0,T)\times \mathbb{R}^d_+)}+\norm{D^2 u}{\Leb{s,q}((0,T)\times \mathbb{R}^d_+)}+\norm{\nabla p}{\Leb{s,q}((0,T)\times \mathbb{R}^d_+)}\\
&\leq N \left(\norm{f}{\Leb{s,q}((0,T)\times \mathbb{R}^d_+)} +\norm{Dg}{\Leb{s,q}((0,T)\times \mathbb{R}^d_+)}+\norm{G}{\Leb{s,q}((0,T)\times \mathbb{R}^d_+)}\right)\\
&\relphantom{=}+N(1+\sqrt{\lambda})\left(\norm{u}{\Leb{s,q}((0,T)\times \mathbb{R}^d_+)}+\norm{p}{\Leb{s,q}((0,T)\times \mathbb{R}^d_+)}+\norm{g}{\Leb{s,q}((0,T)\times \mathbb{R}^d_+)}\right)\\
&\relphantom{=}+\frac{N}{\sqrt{\lambda}}\norm{h}{\Leb{s,q}((0,T)\times \mathbb{R}^d_+)}.
\end{align*}
\end{lemma}
\begin{proof}
We prove this lemma by employing an idea of Agmon \cite{A62} (see e.g. \cite[Lemma 5.5]{K07}). Let $\zeta \in C_0^\infty(\mathbb{R})$ be an odd function with respect to the origin. Define 
\begin{align*}
{u}^{\flat}(t,y,x)=(0,u(t,x)\cos(\sqrt{\lambda}y)\zeta(y)),\quad
{p}^{\flat}(t,y,x)=p(t,x)\cos(\sqrt{\lambda}y)\zeta(y).
\end{align*}
Let 
\begin{align*}
{f}^{\flat}(t,y,x)&=\partial_t{u}^{\flat}-\sum_{i,j=1}^d {a}^{ij} D_{x_i x_j} {u}^{\flat}-D_y^2 {u}^{\flat}-\sum_{i=1}^d {b}^i D_{x_i}{u}^{\flat}-\nabla_{y,x} {p}^{\flat},\\
{g}^{\flat}(t,y,x)&=g(t,x) \zeta(y)\cos(\sqrt{\lambda} y).
\end{align*}
Note that 
\begin{align*}
D_y^2 {u}^{\flat}&=\left(0,u(\zeta''(y)\cos(\sqrt{\lambda}y)-2\sqrt{\lambda} \zeta'(y)\sin(\sqrt{\lambda}y)-\lambda \zeta(y) \cos(\sqrt{\lambda}y))\right),\\
\nabla_{y,x}{p}^{\flat}&=((\zeta'(y) \cos(\sqrt{\lambda}y)-\sqrt{\lambda}\sin(\sqrt{\lambda}y)\zeta(y))p(t,x),\nabla_x p(t,x) \cos(\sqrt{\lambda}y)\zeta(y)).
\end{align*}
If we write $({f}^{\flat})'=(({f}^{\flat})^2,\dots,({f}^{\flat})^{d+1})$, then we have
\begin{align*}
({f}^{\flat})^1 &= (\zeta'(y) \cos(\sqrt{\lambda}y)-\sqrt{\lambda}\sin(\sqrt{\lambda}y)\zeta(y))p(t,x), \\
({f}^{\flat})' &=\zeta(y)\cos(\sqrt{\lambda}y)f+\left(2\sqrt{\lambda} \sin(\sqrt{\lambda}y)\zeta'(y)-\zeta''(y)\cos(\sqrt{\lambda} y)\right)u.
\end{align*}

On the other hand, note that 
\[    \partial_t {g}^{\flat}(t,y,x)= \Div_{y,x} {G}^{\flat}(t,y,x),\]
where 
\[
   {G}^{\flat}(t,y,x)=\left(h(t,x)\int_{-\infty}^y \zeta(z) \cos(\sqrt{\lambda} z)dz,G(t,x)\zeta(y)\cos(\sqrt{\lambda} y) \right).
\]

Since $\zeta$ is odd and has compact support, there exists $\eta_\lambda \in \Leb{\infty}(\mathbb{R})$ that has the same compact support and 
\begin{equation}\label{eq:zeta-expression}
  \int_{-\infty}^y \zeta(z)\cos(\sqrt{\lambda} z)dz=\frac{1}{\sqrt{\lambda}} \eta_\lambda(y),
\end{equation}
where $\norm{\eta_\lambda}{\Leb{\infty}}\leq N$ and the constant $N$ is independent of $\lambda$. 
Since 
\[ \int_{\mathbb{R}} |\zeta(y)\cos(\sqrt{\lambda} y)|^p\myd{y} \geq c>0,\]
where $c$ is independent of $\lambda$, it follows from Theorem \ref{thm:VMO-half} that  
\begin{align*}
&\lambda \norm{u}{\Leb{s,q}((0,T)\times \mathbb{R}^d_+)} + \norm{\partial_t u}{\Leb{s,q}((0,T)\times \mathbb{R}^d_+)}+\norm{D^2 u}{\Leb{s,q}((0,T)\times \mathbb{R}^d_+)}\\
&\leq N\left(\norm{D_y^2 {u}^{\flat}}{\Leb{s,q}((0,T)\times \mathbb{R}\times \mathbb{R}^d_+)}+\norm{\partial_t {u}^{\flat}}{\Leb{s,q}((0,T)\times \mathbb{R}\times \mathbb{R}^d_+)}+\norm{D^2_x {u}^{\flat}}{\Leb{s,q}((0,T)\times \mathbb{R}\times \mathbb{R}^d_+)}\right)\\
&\leq N\left( \norm{{f}^{\flat}}{\Leb{s,q}((0,T)\times \mathbb{R}\times \mathbb{R}^d_+)}+\norm{D_{x,y}{g}^{\flat}}{\Leb{s,q}((0,T)\times \mathbb{R}\times \mathbb{R}^d)} +\norm{{G}^{\flat}}{\Leb{s,q}((0,T)\times \mathbb{R}\times \mathbb{R}^d_+)}\right.\\
&\relphantom{=}\quad\left. +\norm{{u}^{\flat}}{\Leb{s,q}((0,T)\times\mathbb{R}\times\mathbb{R}^d_+)}\right)
\end{align*}
for some constant $N=N(d,s,q,\nu,R_0)>0$. By \eqref{eq:zeta-expression}, we have
\begin{align*}
\norm{{f}^{\flat}}{\Leb{s,q}((0,T)\times \mathbb{R}\times \mathbb{R}^d_+)}&\apprle (1+\sqrt{\lambda})\left(\norm{u}{\Leb{s,q}((0,T)\times \mathbb{R}^d_+)}+\norm{p}{\Leb{s,q}((0,T)\times \mathbb{R}^d_+)}\right)\\
&\relphantom{=}+ \norm{f}{\Leb{s,q}((0,T)\times \mathbb{R}^d_+)},\\
\norm{D_{x,y}{g}^{\flat}}{\Leb{s,q}((0,T)\times \mathbb{R}\times \mathbb{R}^d_+)}&\apprle \norm{Dg}{\Leb{s,q}((0,T)\times \mathbb{R}^d_+)}+(1+\sqrt{\lambda}) \norm{g}{\Leb{s,q}((0,T)\times \mathbb{R}^d_+)},\\
\norm{{G}^{\flat}}{\Leb{s,q}((0,T)\times \mathbb{R}\times \mathbb{R}^d_+)}&\apprle \frac{1}{\sqrt{\lambda}} \norm{h}{\Leb{s,q}((0,T)\times \mathbb{R}^d_+)} + \norm{G}{\Leb{s,q}((0,T)\times \mathbb{R}^d_+)}.
\end{align*}
Hence it follows that 
\begin{align*}
&\lambda \norm{u}{\Leb{s,q}((0,T)\times \mathbb{R}^d_+)} +\norm{\partial_t u}{\Leb{s,q}((0,T)\times \mathbb{R}^d_+)}+\norm{D^2u}{\Leb{s,q}((0,T)\times \mathbb{R}^d_+)}\\
&\leq N \left(\norm{f}{\Leb{s,q}((0,T)\times \mathbb{R}^d_+)} + (1+\sqrt{\lambda})(\norm{u}{\Leb{s,q}((0,T)\times \mathbb{R}^d_+)}+\norm{p}{\Leb{s,q}((0,T)\times \mathbb{R}^d_+)})\right)\\
&\relphantom{=}+ N\left(\norm{Dg}{\Leb{s,q}((0,T)\times \mathbb{R}^d_+)}+ (1+\sqrt{\lambda})\norm{g}{\Leb{s,q}((0,T)\times \mathbb{R}^d_+)}\right)\\
&\relphantom{=}+ N \left(\frac{1}{\sqrt{\lambda}}\norm{h}{\Leb{s,q}((0,T)\times \mathbb{R}^d_+)} + \norm{G}{\Leb{s,q}((0,T)\times \mathbb{R}^d_+)}\right)
\end{align*}
for some constant $N=N(d,s,q,\nu,R_0)>0$. The desired estimate follows by an interpolation inequality. This completes the proof of Lemma \ref{lem:lambda-estimate}.
\end{proof}
  
\begin{proof}[Proof of Theorem \ref{thm:B}]
For simplicity, we write
\[ \norm{v}{X_{s,q;r}}=\norm{\partial_t v}{\Leb{s,q}(Q_r^+)} +\norm{D^2 v}{\Leb{s,q}(Q_r^+)} \]
and
\[ \norm{v}{X_{s,q}}=\norm{\partial_t v}{\Leb{s,q}((-R^2,\infty)\times \mathbb{R}^d_+)} +\norm{D^2 v}{\Leb{s,q}((-R^2,\infty)\times \mathbb{R}^d_+)}. \]

Define $w=\tilde{u}-v$. Then $(w,\tilde{p}) \in \Sob{1,2}{s,q}((-R^2,0)\times \mathbb{R}^d_+)^d\times \Sob{0,1}{s,q}((-R^2,0)\times \mathbb{R}^d_+)$ satisfies 
\begin{equation*}
\left\{
\begin{aligned}
\partial_t w -\mathcal{L} w+\nabla \tilde{p}+\lambda w&=\tilde{f}+(\mathcal{L}-\Delta)v+\lambda\tilde{u}\\
\Div w&=\tilde{g}+\Div((I-\widehat{C})\tilde{u}-v)
\end{aligned}
\right.
\end{equation*}
in $(-R^2,0)\times \mathbb{R}^d_+$, $w(-R^2,\cdot)=0$, and 
\[  w^d=0,\quad D_d w^k=0\quad \text{for all } k=1,2,\dots,d-1\quad \text{on } (-R^2,0]\times \mathbb{R}^{d-1}\times\{0\},\]
where $\tilde{f}$ and $\tilde{g}$ are defined in \eqref{eq:f-tilde-g-tilde}. 
Recall from \eqref{eq:G-hat-defn} that $\partial_t \hat{g}=\Div \widehat{G}$ in $Q_R^+$. Then we have
\begin{align*}
  &\partial_t \left[ \tilde{g} + \Div((I-\widehat{C})\tilde{u}-v)\right]\\
  &=[\hat{g}\zeta + \widehat{C}\tilde{u}\cdot \nabla \zeta]_t + \Div((I-\widehat{C})\tilde{u}_t-v_t)\\
  &=\hat{g}_t \zeta + \hat{g}\zeta_t + \widehat{C}\tilde{u}_t \cdot \nabla \zeta + \widehat{C}\tilde{u} \cdot \nabla \zeta_t + \Div((I-\widehat{C})\tilde{u}_t-v_t)\\
  &=-\hat{G}\cdot\nabla \zeta + \hat{g}\zeta_t + \widehat{C} \tilde{u}_t \cdot \nabla \zeta + \widehat{C}\tilde{u} \cdot \nabla \zeta_t + \Div((I-\widehat{C})\tilde{u}_t-v_t+\widehat{G}\zeta).
\end{align*}
Since we will choose $\lambda$ sufficiently large, we may assume that $\lambda \geq 1$. Then it follows from Lemma \ref{lem:lambda-estimate} that
\begin{align*}
&\lambda\norm{w}{\Leb{s,q}((-R^2,0)\times\mathbb{R}^d_+)}+\norm{\partial_t w}{\Leb{s,q}((-R^2,0)\times\mathbb{R}^d_+)}\\
&\relphantom{=}+\norm{D^2 w}{\Leb{s,q}((-R^2,0)\times\mathbb{R}^d_+)}+\norm{\nabla \hat{p}}{\Leb{s,q}(Q_r^+)}\\
&\apprle \norm{\tilde{f}+(\mathcal{L}-\Delta)v+\lambda\tilde{u}}{\Leb{s,q}((-R^2,0)\times\mathbb{R}^d_+)}+\sqrt{\lambda}\norm{w}{\Leb{s,q}((-R^2,0)\times\mathbb{R}^d_+)}\\
&\relphantom{=}+\sqrt{\lambda}\norm{\hat{p}}{\Leb{s,q}(Q_R^+)}+\norm{D[\tilde{g}+\Div((I-\widehat{C})\tilde{u}-v)]}{\Leb{s,q}((-R^2,0)\times\mathbb{R}^d_+)}\\
&\relphantom{=}+\sqrt{\lambda}\norm{\tilde{g}+\Div((I-\widehat{C})\tilde{u}-v)}{\Leb{s,q}((-R^2,0)\times\mathbb{R}^d_+)}\\
&\relphantom{=}+\frac{1}{\sqrt{\lambda}}\norm{-\widehat{G}\cdot \nabla \zeta +\hat{g}\zeta_t+\widehat{C}\partial_t \tilde{u}\cdot \nabla \zeta+\widehat{C}\tilde{u}\cdot \nabla\zeta_t}{\Leb{s,q}((-R^2,0)\times\mathbb{R}^d_+)}\\
&\relphantom{=}+\norm{(I-\widehat{C})\tilde{u}_t-v_t+\widehat{G}\zeta}{\Leb{s,q}((-R^2,0)\times\mathbb{R}^d_+)}\\
&=\mathrm{I}+\mathrm{II}+\mathrm{III}+\mathrm{IV}+\mathrm{V}+\mathrm{VI}+\mathrm{VII},
\end{align*}
where the implicit constant is independent of $\lambda$, $r$, and $R$. Now we estimate each term. From now on, we may assume that $R\leq 1$. \bigskip

\emph{Estimates on $\mathrm{I}$.} By interpolation inequalities, for each $\varepsilon\in(0,1)$, we have
\begin{align*}
\mathrm{I}&\apprle \norm{\hat{f}}{\Leb{s,q}(Q_R^+)}+\frac{1}{(R-r)^2} \norm{\hat{u}}{\Leb{s,q}(Q_R^+)} + \norm{D\phi}{\Leb{\infty}(B_R')}\norm{\nabla \hat{p}}{\Leb{s,q}(Q_R^+)}\\
&\relphantom{=}+\frac{1}{R-r}\left(\norm{D\hat{u}}{\Leb{s,q}(Q_R^+)}+\norm{\hat{p}}{\Leb{s,q}(Q_R^+)}\right)+\lambda\norm{\hat{u}}{\Leb{s,q}(Q_R^+)}\\
&\relphantom{=}+\norm{v}{X_{s,q}}+\norm{Dv}{\Leb{s,q}((-R^2,\infty)\times\mathbb{R}^d_+)}\\
&\apprle \varepsilon \norm{D^2\hat{u}}{\Leb{s,q}(Q_R^+)} + \norm{D\phi}{\Leb{\infty}(B_R')}\norm{\nabla \hat{p}}{\Leb{s,q}(Q_R^+)}\\
&\relphantom{=}+\norm{v}{X_{s,q}}+\norm{Dv}{\Leb{s,q}((-R^2,\infty)\times\mathbb{R}^d_+)}\\
&\relphantom{=}+\norm{\hat{f}}{\Leb{s,q}(Q_R^+)}+\frac{1}{R-r}\norm{\hat{p}}{\Leb{s,q}(Q_R^+)}+\left(\lambda+ \frac{1}{\varepsilon(R-r)^2}\right)\norm{\hat{u}}{\Leb{s,q}(Q_R^+)}.
\end{align*}

\emph{Estimates on $\mathrm{IV}$.} Since $\tilde{g}=\hat{g}\zeta + (\widehat{C}\hat{u})\cdot \nabla \zeta$, we have
\begin{align*}
\mathrm{IV} &\apprle \norm{D\tilde{g}}{\Leb{s,q}((-R^2,0)\times\mathbb{R}^d_+)}+\norm{D^2 v}{\Leb{s,q}((-R^2,0)\times\mathbb{R}^d_+)}\\
&\relphantom{=}+\norm{D \Div((I-\widehat{C})\tilde{u})}{\Leb{s,q}((-R^2,0)\times\mathbb{R}^d_+)}\\
&\apprle \norm{D\hat{g}}{\Leb{s,q}(Q_R^+)} + \frac{1}{R-r}\norm{\hat{g}}{\Leb{s,q}(Q_R^+)}+\frac{1}{(R-r)^2}\norm{\hat{u}}{\Leb{s,q}(Q_R^+)}+\frac{1}{R-r}\norm{D\hat{u}}{\Leb{s,q}(Q_R^+)} \\
&\relphantom{=}+\norm{v}{X_{s,q}}+\norm{I-\widehat{C}}{\Leb{\infty}(B_R)}\norm{D^2 \hat{u}}{\Leb{s,q}(Q_R^+)}.
\end{align*}
Hence it follows from interpolation inequality that for each $\varepsilon\in(0,1)$, we have
\begin{align*}
\mathrm{IV} &\apprle (\norm{I-\widehat{C}}{\Leb{\infty}(B_R)}+\varepsilon)\norm{D^2 \hat{u}}{\Leb{s,q}(Q_R^+)} +\norm{v}{X_{s,q}}\\
&\relphantom{=}+\norm{D\hat{g}}{\Leb{s,q}(Q_R^+)} + \frac{1}{R-r}\norm{\hat{g}}{\Leb{s,q}(Q_R^+)}+\frac{1}{\varepsilon(R-r)^2}\norm{\hat{u}}{\Leb{s,q}(Q_R^+)}.
\end{align*}


\emph{Estimates on $\mathrm{V}$.}  Since $\tilde{g}=\hat{g}\zeta + (\widehat{C}\hat{u})\cdot \nabla \zeta$, we have 
\begin{align*}
\mathrm{V}&\apprle \sqrt{\lambda}\left(\norm{\tilde{g}}{\Leb{s,q}(Q_R^+)} +  \norm{Dv}{\Leb{s,q}((-R^2,\infty)\times\mathbb{R}^d_+)}+ \norm{D((I-\widehat{C})\tilde{u})}{\Leb{s,q}((-R^2,\infty)\times\mathbb{R}^d_+)}\right)\\
&\apprle \sqrt{\lambda}\left(\norm{\hat{g}}{\Leb{s,q}(Q_R^+)}+ \norm{Dv}{\Leb{s,q}((-R^2,\infty)\times\mathbb{R}^d_+)}+\frac{1}{R-r}\norm{\hat{u}}{\Leb{s,q}(Q_R^+)}\right.\\
&\relphantom{=}\left.+\norm{I-\widehat{C}}{\Leb{\infty}(B_R)}\norm{D\hat{u}}{\Leb{s,q}(Q_R^+)}\right).
\end{align*}
Since $R\leq 1$, it follows from interpolation inequalities that for each $\varepsilon\in(0,1)$, we have 
\begin{align*}
\mathrm{V}&\apprle \varepsilon \norm{D^2\hat{u}}{\Leb{s,q}(Q_R^+)}+\sqrt{\lambda}\left(\norm{\hat{g}}{\Leb{s,q}(Q_R^+)} +  \norm{Dv}{\Leb{s,q}((-R^2,\infty)\times\mathbb{R}^d_+)}\right)\\
&\relphantom{=}+\left(\frac{\sqrt{\lambda}}{R-r}+\frac{\lambda}{\varepsilon} \right)\norm{\hat{u}}{\Leb{s,q}(Q_R^+)}.
\end{align*}

\emph{Estimates on $\mathrm{VI}$.} We have
\begin{align*}
\mathrm{VI}&\apprle  \frac{1}{(R-r)\sqrt{\lambda} } \norm{\partial_t\hat{u}}{\Leb{s,q}(Q_R^+)}+\frac{1}{(R-r)\sqrt{\lambda}}\norm{\widehat{G}}{\Leb{s,q}(Q_R^+)}\\
&\relphantom{=}+\frac{1}{(R-r)^2\sqrt{\lambda} }\norm{\hat{g}}{\Leb{s,q}(Q_R^+)}+\frac{1}{(R-r)^3\sqrt{\lambda}}\norm{\hat{u}}{\Leb{s,q}(Q_R^+)}.
\end{align*}


\emph{Estimates on $\mathrm{VII}$.} We have
\begin{align*}
\mathrm{VII}&\apprle \norm{I-\widehat{C}}{\Leb{\infty}(B_R')}\norm{\partial_t\hat{u}}{\Leb{s,q}(Q_R^+)}\\
&\relphantom{=}+\norm{\partial_t v}{X_{s,q}}+\norm{\hat{G}}{\Leb{s,q}(Q_R^+)}+\frac{1}{(R-r)^2}\norm{\hat{u}}{\Leb{s,q}(Q_R^+)}.
\end{align*}

Choose $\lambda_0>1$ sufficiently large so that $\lambda_0 \geq 2N\sqrt{\lambda_0}$, where $N$ is an implicit constant in the above estimates. For $\lambda \geq \lambda_0$,  we can absorb the term $\mathrm{II}$ into the left-hand side to get
\begin{align*}
\norm{w}{X_{s,q}}+\norm{\nabla \hat{p}}{\Leb{s,q}(Q_r^+)}&\leq \mathrm{I}+ \mathrm{III}+\mathrm{IV}+\mathrm{V}+\mathrm{VI}+\mathrm{VII}.
\end{align*}

Then by \eqref{eq:v-estimate-simplified} and the estimates above, we have
\begin{align}
&\norm{\hat{u}}{X_{s,q;r}}+\norm{\nabla \hat{p}}{\Leb{s,q}(Q_r^+)}\leq \norm{v}{X_{s,q}}+\norm{w}{X_{s,q}}+\norm{\nabla \hat{p}}{\Leb{s,q}(Q_r^+)}\nonumber\\
&\apprle \left(\varepsilon+\norm{I-\hat{C}}{\Leb{\infty}(B_R)}+\frac{1}{(R-r)\sqrt{\lambda}} \right)\norm{\hat{u}}{X_{s,q;R}}+\norm{D\phi}{\Leb{\infty}(B_R')}\norm{\nabla\hat{p}}{\Leb{s,q}(Q_R^+)}\nonumber\\
&\relphantom{=}+\norm{v}{X_{s,q}}+\sqrt{\lambda}\norm{Dv}{\Leb{s,q}((-R^2,\infty)\times\mathbb{R}^d_+)}+\left(\frac{1}{(R-r)\sqrt{\lambda}}+1\right)\norm{\hat{G}}{\Leb{s,q}(Q_R^+)}\nonumber\\
&\relphantom{=}+\left(\frac{\sqrt{\lambda}}{R-r}+\frac{1}{\varepsilon(R-r)^2}+\frac{\lambda}{\varepsilon}+\frac{1}{(R-r)^3\sqrt{\lambda}}\right)\norm{\hat{u}}{\Leb{s,q}(Q_R^+)}\label{eq:u-p-g-estimates}\\
&\relphantom{=}+\left(\frac{1}{R-r}+\frac{1}{(R-r)^2\sqrt{\lambda}}+\sqrt{\lambda} \right)\norm{\hat{g}}{\Leb{s,q}(Q_R^+)}\nonumber\\
&\relphantom{=}+\norm{D\hat{g}}{\Leb{s,q}(Q_R^+)}+\norm{\hat{f}}{\Leb{s,q}(Q_R^+)}+\left(\frac{1}{R-r}+\sqrt{\lambda}\right)\norm{\hat{p}}{\Leb{s,q}(Q_R^+)},\nonumber
\end{align}
where $N$ is independent of $\varepsilon$, $r$, $R$, and $\lambda$. 

Since $r<R$, if we take $\sqrt{\lambda}=(R-r)^{-1}\theta^{-1}$, then it follows from  \eqref{eq:v-estimate-simplified} and \eqref{eq:u-p-g-estimates} that 
\begin{align*}
&\norm{\hat{u}}{X_{s,q;r}}+\norm{\nabla \hat{p}}{\Leb{s,q}(Q_r^+)}\\
&\leq N\left(\varepsilon+\norm{D\phi}{\Leb{\infty}(B_R')}+\theta+\norm{I-\widehat{C}}{\Leb{\infty}(B_R)} \right)\norm{\hat{u}}{X_{s,q;R}}\\
&\relphantom{=}+N\norm{D\phi}{\Leb{\infty}(B_R')}\norm{\nabla\hat{p}}{\Leb{s,q}(Q_R^+)} +N_{\theta,\varepsilon}\left(\frac{1}{(R-r)^2}+\frac{1}{(R-r)^{2+b}}\right)\norm{\hat{u}}{\Leb{s,q}(Q_R^+)}\\
&\relphantom{=}+\frac{N}{R-r}\left(\norm{\hat{p}}{\Leb{s,q}(Q_R^+)}+\norm{\hat{g}}{\Leb{s,q}(Q_R^+)}\right)+N\norm{D\hat{g}}{\Leb{s,q}(Q_R^+)}\\
&\relphantom{=}+N\norm{\hat{f}}{\Leb{s,q}(Q_R^+)}+N\norm{\hat{G}}{\Leb{s,q}(Q_R^+)},
\end{align*}
where $N$ is independent of $\varepsilon$, $r$, $R$, $\lambda$, and $N_{\theta,\varepsilon}$ is a constant that depends on $\theta$ and $\varepsilon$ but independent of $r$, $R$, and $\lambda$.

Define
\[  r_k= R-\frac{R-r}{2^{k}} ,\quad \mathsf{A}_k=\norm{\hat{u}}{X_{s,q;r_k}},\quad \mathsf{B}_k=\norm{\nabla\hat{p}}{\Leb{s,q}(Q_{r_k}^+)},\quad\text{and}\quad\mathsf{C}_k=\norm{\hat{u}}{\Leb{s,q}(Q_{r_k}^+)}.\]
If we choose $R^*\leq 1$, $\varepsilon$, $\theta$ sufficiently small so that 
\[ N(\varepsilon +\theta +\norm{I-\widehat{C}}{\Leb{\infty}(B_{R^*})} +\norm{D\phi}{\Leb{\infty}(B_{R^*}')})\leq \frac{1}{3^{2+b}},\]
then  for $0<r<R\leq R^*$, we have 
\begin{equation}\label{eq:final-iteration}
\begin{aligned}
\mathsf{A}_k +\mathsf{B}_k &\leq \frac{1}{3^{2+b}}(\mathsf{A}_{k+1}+\mathsf{B}_{k+1}) +N\frac{2^{(2+b) k}}{(R-r)^{2+b}}\mathsf{C}_{k+1}\\
&\relphantom{=}+N\frac{ 2^k}{R-r}\left(\norm{\hat{p}}{\Leb{s,q}(Q_R^+)}+\norm{\hat{g}}{\Leb{s,q}(Q_R^+)}\right)+N\norm{D\hat{g}}{\Leb{s,q}(Q_R^+)}\\
&\relphantom{=}+N\norm{\hat{f}}{\Leb{s,q}(Q_R^+)}+N\norm{\hat{G}}{\Leb{s,q}(Q_R^+)}.
\end{aligned}
\end{equation}
Multiplying both sides of the inequality \eqref{eq:final-iteration} by $3^{-(2+b)k}$ and taking summation over $k=0,1,\dots,$ we get 
\begin{align*}
\sum_{k=0}^\infty \left(\frac{1}{3^{2+b}}\right)^k (\mathsf{A}_k+\mathsf{B}_k)&\leq \sum_{k=0}^\infty \left(\frac{1}{3^{2+b}}\right)^{k+1} (\mathsf{A}_{k+1}+\mathsf{B}_{k+1})\\
&\relphantom{=}+\frac{N}{(R-r)^{2+b}}\norm{\hat{u}}{\Leb{s,q}(Q_R^+)}+N\norm{\hat{f}}{\Leb{s,q}(Q_R^+)}\\
&\relphantom{=}+N \norm{\hat{G}}{\Leb{s,q}(Q_R^+)}+\frac{N}{R-r}\left(\norm{\hat{p}}{\Leb{s,q}(Q_R^+)}+\norm{\hat{g}}{\Sob{0,1}{s,q}(Q_R^+)}\right).
\end{align*}
Hence by cancelling the summation on the right-hand side, we get 
\begin{align*}
&\norm{\partial_t \hat{u}}{\Leb{s,q}(Q_r^+)}+\norm{D^2\hat{u}}{\Leb{s,q}(Q_r^+)}+\norm{\nabla \hat{p}}{\Leb{s,q}(Q_r^+)}\\
&\leq N(d,s,q,r,R,M,\mathcal{A})\left(\norm{\hat{u}}{\Leb{s,q}(Q_R^+)}+\norm{\hat{p}}{\Leb{s,q}(Q_R^+)}+\norm{\hat{f}}{\Leb{s,q}(Q_R^+)}\right.\\
&\relphantom{=}\left.+\norm{\hat{g}}{\Sob{0,1}{s,q}(Q_R^+)}+\norm{\hat{G}}{\Leb{s,q}(Q_R^+)}\right).
\end{align*}
Then the desired result follows by a standard covering argument. This completes the proof of Theorem \ref{thm:B}.
\end{proof}



\subsection*{Acknowledgements}
The second author thanks Aidan Backus for the discussion related to Appendix \ref{sec:principal}.
\appendix
\section{Navier boundary condition}\label{sec:principal}

Let $n$ be the outward unit normal vector to $\partial\Omega$. For a vector $\tau$ in $\mathbb{R}^d$, let $D_\tau$ denote the directional derivative with respect to $\tau$. If $\tau$ is tangent to $\partial\Omega$, let $\nabla_\tau$ denote the covariant directional derivative with respect to $\tau$, computed using the Riemannian metric on $\partial\Omega$. Then by Gauss-Codazzi theorem, there exists a symmetric $2$-tensor field $\mathcal{A}_w$ on $\partial\Omega$ such that for every vector field $v$ tangent to $\partial\Omega$, and every unit tangent vector $\tau$ to $\partial\Omega$, we have
\begin{equation}\label{eq:Weingarten}
 (D_\tau v)_i = (\nabla_\tau v)_i -(\tau \cdot \mathcal{A}_wv)n_i.
\end{equation}
The tensor $\mathcal{A}_w$ is called the Weingarten map or shape operator of $\partial\Omega$. 

\begin{lemma}
Let $\Omega$ be a $C^{2,1}$-domain in $\mathbb{R}^d$. If $u\cdot n=0$ on $\partial\Omega$, then 
\begin{equation}\label{eq:generalized-Navier}
 [\mathbb{D}(u)n+\mathcal{A}_w u]_{\mathrm{tan}}=\frac{1}{2}\bm{\omega}n.
\end{equation}
\end{lemma}
\begin{proof}
Since $u\cdot n=0$ on $\partial\Omega$, $u$ is tangent to $\partial\Omega$. For any vector field $\tau$ which is tangent to $\partial\Omega$, we have 
\[ 2\tau \cdot (\mathbb{D}(u)n)=(\bm{\omega}n)\cdot\tau + 2(D_\tau u)\cdot n. \]
By \eqref{eq:Weingarten}, we have 
\[ (D_\tau u)\cdot n = (\nabla_\tau u)\cdot n -(\tau \cdot \mathcal{A}_w u)=-\tau \cdot \mathcal{A}_w u. \]
This implies \eqref{eq:generalized-Navier}.
\end{proof}

\bibliographystyle{amsplain}

\begin{thebibliography}{10}

\bibitem{A62}
S.~Agmon, \emph{On the eigenfunctions and on the eigenvalues of general
  elliptic boundary value problems}, Comm. Pure Appl. Math. \textbf{15} (1962),
  119--147. 

\bibitem{BB99}
J.-L. Barrat and L.~Bocquet, \emph{Large slip effect at a {N}onwetting
  fluid-solid interface}, Physical Review Letters \textbf{82} (1999), no.~23,
  4671--4674.

\bibitem{BGL16}
A.~Bonito, J.-L. Guermond, and S.~Lee, \emph{Numerical simulations of bouncing
  jets}, Internat. J. Numer. Methods Fluids \textbf{80} (2016), no.~1, 53--75.

\bibitem{CK20}
T.~Chang and K.~Kang, \emph{On {C}accioppoli's inequalities of {S}tokes
  equations and {N}avier-{S}tokes equations near boundary}, J. Differential
  Equations \textbf{269} (2020), no.~9, 6732--6757. 

\bibitem{CSTY09}
C.-C. Chen, R.~M. Strain, T.-P. Tsai, and H.-T. Yau, \emph{Lower bounds on the
  blow-up rate of the axisymmetric {N}avier-{S}tokes equations. {II}}, Comm.
  Partial Differential Equations \textbf{34} (2009), no.~1-3, 203--232.

\bibitem{CLT23}
H.~Chen, S.~Liang, and T.-P. Tsai, \emph{Gradient estimates for the
  non-stationary {S}tokes system with the {N}avier boundary condition}, Comm.
  Pure Appl. Anal. \textbf{23} (2024), no.~10, 1483--1505.

\bibitem{CLT24}
\bysame, \emph{Poisson kernel and blow-up of the second derivatives near the
  boundary for {S}tokes equations with {N}avier boundary condition},
  arXiv:2406.15995, 2024.

\bibitem{DHP07}
R.~Denk, M.~Hieber, and J.~Pr\"{u}ss, \emph{Optimal {$L^p$}-{$L^q$}-estimates
  for parabolic boundary value problems with inhomogeneous data}, Math. Z.
  \textbf{257} (2007), no.~1, 193--224. 

\bibitem{DM04}
M.~Dindo\v{s} and M.~Mitrea, \emph{The stationary {N}avier-{S}tokes system in
  nonsmooth manifolds: the {P}oisson problem in {L}ipschitz and {$C^1$}
  domains}, Arch. Ration. Mech. Anal. \textbf{174} (2004), no.~1, 1--47.

\bibitem{DKP22}
H.~Dong, D.~Kim, and T.~Phan, \emph{Boundary {L}ebesgue mixed-norm estimates
  for non-stationary {S}tokes systems with {VMO} coefficients}, Comm. Partial
  Differential Equations \textbf{47} (2022), no.~8, 1700--1731. 

\bibitem{DKr19}
H.~Dong and N.~V. Krylov, \emph{Fully nonlinear elliptic and parabolic
  equations in weighted and mixed-norm {S}obolev spaces}, Calc. Var. Partial
  Differential Equations \textbf{58} (2019), no.~4, Paper No. 145, 32.

\bibitem{DK23}
H.~Dong and H.~Kwon, \emph{Interior and boundary mixed norm derivative
  estimates for nonstationary stokes equations}, arXiv:2308.09220, 2023.

\bibitem{DP21}
H.~Dong and T.~Phan, \emph{Mixed-norm {$L_p$}-estimates for non-stationary
  {S}tokes systems with singular {VMO} coefficients and applications}, J.
  Differential Equations \textbf{276} (2021), 342--367. 

\bibitem{DL22}
R.~Dong and D.~Li, \emph{Interior {$L_p$} regularity for {S}tokes systems},
  arXiv:2212.13698, 2022.


\bibitem{GK12}
G.-M. Gie and J.~P. Kelliher, \emph{Boundary layer analysis of the
  {N}avier-{S}tokes equations with generalized {N}avier boundary conditions},
  J. Differential Equations \textbf{253} (2012), no.~6, 1862--1892.

\bibitem{G18}
G.~Grubb, \emph{Regularity in {$L_p$} {S}obolev spaces of solutions to
  fractional heat equations}, J. Funct. Anal. \textbf{274} (2018), no.~9,
  2634--2660. 

\bibitem{GHJO20}
Y.~Guo, H.~J. Hwang, J.~W. Jang, and Z.~Ouyang, \emph{The {L}andau equation
  with the specular reflection boundary condition}, Arch. Ration. Mech. Anal.
  \textbf{236} (2020), no.~3, 1389--1454. 

\bibitem{HLW14}
F.~Hu, D.~Li, and L.~Wang, \emph{A new proof of {$L^p$} estimates of {S}tokes
  equations}, J. Math. Anal. Appl. \textbf{420} (2014), no.~2, 1251--1264.

\bibitem{HL22}
F.~Hummel and N.~Lindemulder, \emph{Elliptic and parabolic boundary value
  problems in weighted function spaces}, Potential Anal. \textbf{57} (2022),
  no.~4, 601--669. 

\bibitem{J13}
B.~J. Jin, \emph{On the {C}accioppoli inequality of the unsteady {S}tokes
  system}, Int. J. Numer. Anal. Model. Ser. B \textbf{4} (2013), no.~3,
  215--223. 

\bibitem{K06}
J.~P. Kelliher, \emph{Navier-{S}tokes equations with {N}avier boundary
  conditions for a bounded domain in the plane}, SIAM J. Math. Anal.
  \textbf{38} (2006), no.~1, 210--232. 

\bibitem{K07}
N.~V. Krylov, \emph{Parabolic and elliptic equations with {VMO} coefficients},
  Comm. Partial Differential Equations \textbf{32} (2007), no.~1-3, 453--475.

\bibitem{LS03}
E.~Lauga and H.~A. Stone, \emph{Effective slip in pressure-driven stokes flow},
  Journal of Fluid Mechanics \textbf{489} (2003), 55--77.

\bibitem{MV15}
M.~Meyries and M.~Veraar, \emph{Pointwise multiplication on vector-valued
  function spaces with power weights}, J. Fourier Anal. Appl. \textbf{21}
  (2015), no.~1, 95--136. 

\bibitem{MT01}
M.~Mitrea and M.~Taylor, \emph{Navier-{S}tokes equations on {L}ipschitz domains
  in {R}iemannian manifolds}, Math. Ann. \textbf{321} (2001), no.~4, 955--987.

\bibitem{NEBJ05}
C.~Neto, D.~R. Evans, E.~Bonaccurso, H.~J. Butt, and V.~S.~J. Craig,
  \emph{Boundary slip in newtonian liquids: a review of experimental studies},
  Reports on Progress in Physics \textbf{68} (2005), no.~12, 2859.

\bibitem{ST87}
H.-J. Schmeisser and H.~Triebel, \emph{Topics in {F}ourier analysis and
  function spaces}, A Wiley-Interscience Publication, John Wiley \& Sons, Ltd.,
  Chichester, 1987. 

\bibitem{S00}
G.~A. Seregin, \emph{Some estimates near the boundary for solutions to the
  non-stationary linearized {N}avier-{S}tokes equations}, Zap. Nauchn. Sem.
  S.-Peterburg. Otdel. Mat. Inst. Steklov. (POMI) \textbf{271} (2000),
  204--223, 317. 

\bibitem{SS14}
G.~A. Seregin and T.~N. Shilkin, \emph{The local regularity theory for the
  {N}avier-{S}tokes equations near the boundary}, Proceedings of the {S}t.
  {P}etersburg {M}athematical {S}ociety. {V}ol. {XV}. {A}dvances in
  mathematical analysis of partial differential equations, Amer. Math. Soc.
  Transl. Ser. 2, vol. 232, Amer. Math. Soc., Providence, RI, 2014,
  pp.~219--244. 

\bibitem{TT97}
P.~A. Thompson and S.~M. Troian, \emph{A general boundary condition for liquid
  flow at solid surfaces}, Nature \textbf{389} (1997), 360--362.

\bibitem{VS13}
V.~Vyalov and T.~Shilkin, \emph{Estimates of solutions to the perturbed
  {S}tokes system}, Zap. Nauchn. Sem. S.-Peterburg. Otdel. Mat. Inst. Steklov.
  (POMI) \textbf{410} (2013), 5--24, 187. 

\bibitem{W15}
J.~Wolf, \emph{On the local regularity of suitable weak solutions to the
  generalized {N}avier-{S}tokes equations}, Ann. Univ. Ferrara Sez. VII Sci.
  Mat. \textbf{61} (2015), no.~1, 149--171. 

\bibitem{XX07}
Y.~Xiao and Z.~Xin, \emph{On the vanishing viscosity limit for the 3{D}
  {N}avier-{S}tokes equations with a slip boundary condition}, Comm. Pure Appl.
  Math. \textbf{60} (2007), no.~7, 1027--1055. 

\end{thebibliography}

\providecommand{\bysame}{\leavevmode\hbox to3em{\hrulefill}\thinspace}
\providecommand{\MR}{\relax\ifhmode\unskip\space\fi MR }
\providecommand{\MRhref}[2]{%
  \href{http://www.ams.org/mathscinet-getitem?mr=#1}{#2}
}
\providecommand{\href}[2]{#2}

\end{document}